\documentclass[11pt,letterpaper]{amsart}

\usepackage{times}
\usepackage[T1]{fontenc}
\usepackage{mathrsfs}
\usepackage{latexsym}
\usepackage[dvips]{graphics}
\usepackage{epsfig}
\usepackage{amsmath,amsfonts,amsthm,amssymb,amscd}
\usepackage{color}
\usepackage[all]{xypic}

\newcommand\be{\begin{equation}}
\newcommand\ee{\end{equation}}
\newcommand\bea{\begin{eqnarray}}
\newcommand\eea{\end{eqnarray}}
\newcommand\bi{\begin{itemize}}
\newcommand\ei{\end{itemize}}
\newcommand\ben{\begin{enumerate}}
\newcommand\een{\end{enumerate}}
\newcommand\bc{\begin{center}}
\newcommand\ec{\end{center}}
\newcommand\ba{\begin{array}}
\newcommand\ea{\end{array}}

\newtheorem{thm}{Theorem}[section]

\newtheorem{cor}[thm]{Corollary}
\newtheorem{lem}[thm]{Lemma}
\newtheorem{prop}[thm]{Proposition}

\newtheorem{defi}[thm]{Definition}

\theoremstyle{remark}
\newtheorem{rek}[thm]{Remark}

\newcommand{\C}{\ensuremath{\mathbb{C}}}
\newcommand{\Z}{\ensuremath{\mathbb{Z}}}

\DeclareMathOperator{\End}{End}
\DeclareMathOperator{\Hom}{Hom}
\DeclareMathOperator{\ad}{ad}
\DeclareMathOperator{\YMH}{YMH}

\DeclareMathOperator{\tr}{tr}
\DeclareMathOperator{\pr}{pr}

\DeclareMathOperator{\rank}{rank}

\DeclareMathOperator{\Lie}{Lie}

\DeclareMathOperator{\im}{im}
\DeclareMathOperator{\SL}{SL}

\DeclareMathOperator{\U}{U}

\newcommand{\starbar}{\bar{*}}

\newcommand{\mhiggs}{\mathcal{M}^{Higgs}}
\newcommand{\totalspace}{\mathcal{A} \times \Omega^0(K \otimes \End(E))}

\newcommand{\G}{\mathcal{G}}

\newcommand{\B}{\mathcal{B}}
\newcommand{\hyperquotient}{ / \negthinspace\negthinspace / \negthinspace\negthinspace /}

\newcommand{\doubleslash}{\bigr/ \negthinspace\negthinspace \bigr/}

\begin{document}

\title{Morse Theory and Hyperk\"ahler Kirwan Surjectivity for Higgs bundles}

\date{\today}

\author[Daskalopoulos]{G.D.~Daskalopoulos}
\address{Department of Mathematics, Brown University, Providence, RI 02912 USA}
\email{daskal@math.brown.edu} 
\thanks{G.D. supported in part by NSF grant DMS-0604930}

\author[Weitsman]{J.~Weitsman}
\address{Department of Mathematics, Northeastern University, Boston, MA 02115}
\email{j.weitsman@neu.edu}
\thanks{J.W. supported in part by NSF grant DMS-0405670 and DMS-0907110}

\author[Wentworth]{R.A.~Wentworth}
\address{Department of Mathematics,
   University of Maryland,
   College Park, MD 20742}
\email{raw@umd.edu}
\thanks{R.W. supported in part by NSF grant DMS-0805797}

\author[Wilkin]{G.~Wilkin}
\address{Department of Mathematics, Johns Hopkins University, Baltimore, MD 21218}
\email{graeme@math.jhu.edu}

\date{\today}

\subjclass[2000]{ Primary: 53C26 ; Secondary: 53D20 }

\begin{abstract}
This paper uses Morse-theoretic techniques to compute the equivariant Betti numbers of the space of semistable rank two degree zero Higgs bundles over a compact Riemann surface, a method in the spirit of Atiyah and Bott's original approach for semistable holomorphic bundles. This leads to a natural proof that the hyperk\"ahler Kirwan map is surjective for the non-fixed determinant case.
\end{abstract}

\maketitle

\tableofcontents

\section{Introduction}

The moduli space of semistable holomorphic bundles over a compact Riemann surface is a well-studied object in algebraic geometry. The seminal paper of Atiyah and Bott introduced a new method for computing the cohomology of this space: The equivariant Morse theory of the Yang-Mills functional.  This  and subsequent work  provides substantial information on  its cohomology ring. Also of interest is the \emph{moduli space of semistable Higgs bundles}.  The purpose of this paper is to develop an equivariant Morse theory on the (singular) space of Higgs bundles in order to carry out the Atiyah and Bott program for the case of rank $2$.

The precise setup is as follows. Let $E$ be a complex Hermitian vector bundle of rank $n$ and degree $d_E$ over a compact Riemann surface $M$ of genus $g$.   Let $\mathcal{A}(2, d_E)$ denote the space of Hermitian connections on $E$, and $\mathcal{A}_0(2, d_E)$ the space of traceless Hermitian connections (which can be identified with the space of holomorphic structures on $E$ without or with a fixed determinant bundle). We use $\End(E)$ to denote the bundle of endomorphisms of $E$, $\End_0(E)$  the subbundle of trace-free endomorphisms, and $\ad(E) \subset \End(E)$ (resp.\ $\ad_0(E) \subset \End_0(E)$)  the subbundle of endomorphisms that are skew adjoint with respect to the Hermitian metric.

Let 
$$\mathcal{B}(2, d_E) = \{ (A, \Phi) \in \mathcal{A}(2, d_E) \times \Omega^0(\End(E) \otimes K) : d_A'' \Phi = 0 \}$$
be the space of \emph{Higgs bundles} of degree $d_E$ and rank $n$ over $M$  and let 
$$\mathcal{B}_0(2, d_E) = \{ (A, \Phi) \in \mathcal{A}_0(2, d_E) \times \Omega^0(\End_0(E) \otimes K) : d_A'' \Phi = 0 \}$$
 denote the space of  Higgs bundles with \emph{fixed determinant}. Let $\mathcal{G}$ (resp. $\mathcal{G}^\C$) denote the gauge group of $E$ with structure group $U(2)$ (resp. $GL(2)$) for the non-fixed determinant case, and  $\mathcal{G}_0$ (resp. $\mathcal{G}_0^\C$) the gauge groups with structure group $SU(2)$ (resp. $SL(2)$) for the fixed determinant case.  The  action of these groups on the space of Higgs bundles is given by
\begin{equation}
g \cdot (A, \Phi) = ( g^{-1} A'' g + g^* A' (g^*)^{-1}+ g^{-1} d'' g - (d' g^*) (g^*)^{-1}, g^{-1} \Phi g) ,
\end{equation}
where $A''$ and $A'$ denote the $(0,1)$ and $(1,0)$ parts of the connection form $A$.

The cotangent bundle $pr: T^* \mathcal{A}(2, d_E)\to \mathcal{A}(2, d_E)$
is naturally   
$$T^* \mathcal{A}(2, d_E)\simeq \mathcal{A}(2, d_E) \times \Omega^0(\End(E) \otimes K)$$
 and this gives rise to a  hyperk\"ahler structure preserved by the action of $\mathcal{G}$ (cf.\ \cite{Hitchin87}). The moment maps for this action are 
\begin{align*} \label{eqn:moment-maps}
\begin{split}
\mu_1 & = F_A + [ \Phi, \Phi^*] \\
\mu_2 & = -i \left( d_A'' \Phi + d_A' \Phi^* \right) \\
\mu_3 & = -d_A'' \Phi + d_A' \Phi^* 
\end{split}
\end{align*}
In the sequel, we refer to $\mu_{\mathbb C}=\mu_2+i\mu_3=2i d_A^{\prime\prime}\Phi$ as the \emph{complex moment map}.
The hyperk\"ahler quotient $T^* \mathcal{A}(2, d_E) \hyperquotient \mathcal{G}$ is the space
\begin{equation*}
T^* \mathcal{A}(2, d_E) \hyperquotient \mathcal{G} : = \mu_1^{-1} (\alpha) \cap \mu_2^{-1}(0) \cap \mu_3^{-1}(0) / \mathcal{G} ,
\end{equation*}
where $\alpha$ is a constant multiple of the identity (depending on $d_E$) chosen so that $\mu_1 = \alpha$ minimizes the \emph{Yang-Mills-Higgs} functional 
$$\YMH(A,\Phi) = \| F_A + [\Phi, \Phi^*] \|^2$$

In the following $\mathcal{B}$ or  $\mathcal{B}_0$ (resp.\ $\mathcal A$ or ${\mathcal A}_0$) will often be used to denote the space of Higgs bundles (resp.\ connections) with non-fixed or fixed determinant, and the extra notation will be omitted if the meaning is clear from the context. Let $\mathcal{B}^{st}$ (resp. $\mathcal{B}^{ss}$) denote the space of \emph{stable} (resp. \emph{semistable}) Higgs bundles, those for which every $\Phi$-invariant holomorphic subbundle $F \subset E$ satisfies 
$$\frac{\deg(F)}{\rank(F)} < \frac{\deg(E)}{\rank(E)}\quad \left(\text{resp. } \frac{\deg(F)}{\rank(F)} \leq \frac{\deg(E)}{\rank(E)}\right)
$$
Similarly for $\mathcal{B}_0^{st}$ and ${\mathcal{B}_0}^{ss}$.
Let $\left< u, v \right> = \int_M tr\{u\starbar v\}$ be the $L^2$ inner product on $\Omega^0(\ad(E))$, with associated norm $\left\| u \right\|^2 = \left< u, u \right>$.  The functional  $\YMH$ is defined on $\mathcal{B}$ and $\mathcal{B}_0$, and $\mu^{-1}(\alpha) \cap \mu_\C^{-1}(0)$ is the subset of Higgs bundles that minimize $\YMH$.

Theorems of Hitchin  \cite{Hitchin87} and Simpson  \cite{Simpson88} identify the hyperk\"ahler quotient 
$$\left\{\B_{min}=\mu_1^{-1}(\alpha) \cap \mu_\C^{-1}(0) \right\}/ \mathcal{G}$$
 with the \emph{moduli space of semistable Higgs bundles} of rank $n$, degree $d_E$ and non-fixed determinant, $\mhiggs(2, d_E) = \mathcal{B}^{ss} \doubleslash \mathcal{G}^\C$, and similarly in the fixed determinant case $\mhiggs_0(2,d_E) = \mathcal{B}_0^{ss} \doubleslash {\mathcal{G}_0}^\C$.
Since $-2i d_A'' \Phi = \mu_2 + i \mu_3$, this hyperk\"ahler quotient can be viewed as a symplectic quotient of the singular space of Higgs bundles
\begin{equation*}
T^* \mathcal{A} \hyperquotient \mathcal{G} = \left( \mathcal{B} \cap \mu_1^{-1}(\alpha) \right) / \mathcal{G} 
\end{equation*}

This paper uses the equivariant Morse theory of the functional $\YMH$ on the space $\mathcal{B}$ and $\mathcal{B}_0$  to study the topology of the moduli space of rank $2$ Higgs bundles for both fixed and non-fixed determinant and both degree zero and odd degree. The main results are the following.
\begin{thm}  \label{thm:bettinumbers1}
For the degree zero case, we have the following formulae for the equivariant Poincar\'e polynomials. For the fixed determinant case,
\begin{align}
\begin{split}
P_t^{\mathcal{G}} (\mathcal{B}_0^{ss}(2, 0)) = & P_t(B\mathcal{G}) - \sum_{d=1}^\infty t^{2\mu_d} \frac{(1+t)^{2g}}{1-t^2} \label{eqn:equivariantbetti-fixed0} \\ 
 & + \sum_{d=1}^{g-1} t^{2\mu_d} P_t(\widetilde{S}^{2g-2d-2} M),
 \end{split}
\end{align}
and for the non fixed determinant case,
\begin{align}
\begin{split}
P_t^{\mathcal{G}}(\mathcal{B}^{ss}(2, 0)) = & P_t(B\mathcal{G}) - \sum_{d=1}^\infty t^{2 \mu_d} \frac{(1+t)^{4g}}{(1-t^2)^2} \label{eqn:equivariantbetti-nonfixed0} \\
 & + \sum_{d=1}^{g-1} t^{2\mu_d} P_t(S^{2g-2d-2} M ) \frac{(1+t)^{2g}}{1-t^2},
\end{split}
\end{align}
where $\mu_d = g+2d-1$ and $\widetilde{S}^n M$ denotes the $2^{2g}$-fold cover of the symmetric product $S^n M$ as described in \cite[Sect.\ 7]{Hitchin87}.
\end{thm}

\begin{cor}\label{thm:maintheoremdegreezero}
The equivariant Poincar\'e polynomial of the space of semistable Higgs bundles of rank $2$ and degree zero with fixed determinant over a compact Riemann surface $M$ of genus $g$ is given by
\begin{align*}
P_t^\mathcal{G}(\mathcal{B}_0^{ss}(2, 0)) = & \frac{(1+t^3)^{2g} - (1+t)^{2g} t^{2g+2}}{(1-t^2)(1-t^4)} \nonumber \\
 & \, - t^{4g-4} + \frac{t^{2g+2} (1+t)^{2g}}{(1-t^2)(1-t^4)} + \frac{(1-t)^{2g} t^{4g-4}}{4(1+t^2)} \\
 & \, + \frac{(1+t)^{2g} t^{4g-4}}{2(1 - t^2)} \left( \frac{2g}{t+1} + \frac{1}{t^2 - 1} - \frac{1}{2} + (3-2g) \right) \nonumber \\
 & \, +\frac{1}{2} (2^{2g}-1) t^{4g-4} \left( (1+t)^{2g-2} + (1-t)^{2g-2} - 2 \right) \nonumber
\end{align*}
and in the non-fixed determinant case,
\begin{align*}
P_t^\mathcal{G}(\mathcal{B}^{ss}(2, 0)) = & \frac{(1+t)^{2g}}{(1-t^2)^2 (1-t^4)} \left( (1+t^3)^{2g} - (1+t)^{2g} t^{2g+2} \right) \\
 & \, +\frac{(1+t)^{2g}}{1-t^2} \left( -t^{4g-4} + \frac{t^{2g+2} (1+t)^{2g}}{(1-t^2)(1-t^4)} + \frac{(1-t)^{2g} t^{4g-4}}{4(1+t^2)} \right) \nonumber \\
 & \, + \frac{(1+t)^{4g} t^{4g-4}}{2(1 - t^2)^2} \left( \frac{2g}{t+1} + \frac{1}{t^2 - 1} - \frac{1}{2} + (3-2g) \right) \nonumber
\end{align*}
\end{cor}

The odd degree case was studied by Hitchin  \cite{Hitchin87} using the Morse theory of the functional $ \left\| \Phi \right\|^2$ which appears as (twice) the moment map associated to the $S^1$ action $e^{it}\cdot (A, \Phi) = (A, e^{it} \Phi)$ on the moduli space $\mhiggs_0(2, 1)$. The methods developed in this paper give a new proof of Hitchin's result.
\begin{thm}[cf.\ {\cite[Sect.\  7]{Hitchin87}}]\label{thm:maintheoremdegreeone}
\begin{equation*}
P_t(\mhiggs_0(2, 1)) = P_t(B\mathcal{G}) - \sum_{d=1}^\infty t^{2\mu_d} \frac{(1+t)^{2g}}{1-t^2} + \sum_{d=1}^{g-1} t^{2\mu_d} P_t(\widetilde{S}^{2g-2d-1} M)
\end{equation*}
where $\widetilde{S}^n M$ denotes the $2^{2g}$-fold cover of the symmetric product $S^n M$ as described in \cite[Sect.\ 7]{Hitchin87}. In the non-fixed determinant case,
\begin{multline*}
P_t(\mhiggs(2, 1)) = (1-t^2)P_t(B\mathcal{G}) - \sum_{d=1}^\infty t^{2 \mu_d} (1+t)^{4g} \frac{1}{1-t^2} \\
 + \sum_{d=1}^{g-1} t^{2\mu_d} P_t(S^{2g-2d-1} M \times J_d(M))
\end{multline*}
where $\mu_d = g + 2d - 2$.
\end{thm}

As mentioned above, the moduli space $\mhiggs$ is the hyperk\"ahler quotient of $T^\ast{\mathcal A}$ by the action of $\mathcal{G}$, with associated \emph{hyperk\"ahler Kirwan map}:
\begin{equation*}
\kappa_H : H_\mathcal{G}^*(\totalspace) \rightarrow H_\mathcal{G}^*(\mu_1^{-1}(0) \cap \mu_\C^{-1}(0))
\end{equation*}
induced by the inclusion $\mu_1^{-1}(0) \cap \mu_\C^{-1}(0) \hookrightarrow \totalspace$.
The Morse theory techniques used to prove Theorems \ref{thm:maintheoremdegreezero} and \ref{thm:maintheoremdegreeone} also lead to a natural proof of the following
\begin{thm}
The hyperk\"ahler Kirwan map is surjective for the space of rank $2$ Higgs bundles of non-fixed determinant, for both degree zero and for odd degree.
\end{thm}

For the case of odd degree, surjectivity was previously shown by Hausel and Thaddeus  \cite{HauselThaddeus04} using different methods. The result proved here applies as well to the heretofore unknown degree zero case, and the proof follows naturally from the Morse theory approach used in this paper. In the fixed determinant case, Hitchin's calculation of $P_t(\mhiggs_0(2, 1))$ for a compact genus $2$ surface shows that $b_5(\mhiggs_0(2, 1)) = 34$, however for genus $2$, $b_5(B \mathcal{G}^{SU(2)}) = 4$, hence surjectivity cannot hold in this case.

The most important technical ingredient of this paper is the result of \cite{Wilkin08}  that the gradient flow of $\YMH$ on the spaces $\mathcal{B}$ and $\mathcal{B}_0$ converges to a critical point that corresponds to the graded object of the Harder-Narasimhan-Seshadri filtration of the initial conditions to the gradient flow. The functional $\YMH$  then provides a gauge group equivariant stratification of the spaces $\mathcal{B}$, $\mathcal{B}_0$, and there is a well-defined deformation retraction of each stratum onto an associated set of critical points. This convergence result is sufficient to develop a Morse-type theory on the singular spaces $\mathcal{B}$ and $\mathcal{B}_0$ and to compute the cohomology of the semistable stratums $\mathcal{B}^{ss}$ and $\mathcal{B}_0^{ss}$. It is therefore a consequence of our methods that the lack of Kirwan surjectivity in the fixed determinant case is {\it{not}} due to analytic problems, as one might initially suspect.

More precisely, the results of \cite{Wilkin08}  show that this Morse stratification is the same as the stratification by the type of the Harder-Narasimhan filtration (cf.\ \cite{HauselThaddeus04}). In the case where $\rank(E) = 2$  the strata are enumerated as follows. Given an unstable Higgs pair $(A, \Phi)$, there exists a destabilizing $\Phi$-invariant line bundle $L \subset E$. The quotient $E / L$ is a line bundle (and hence stable), therefore the Harder-Narasimhan filtration is $0 \subset L \subset E$. In this case the type of the Harder-Narasimhan filtration is determined by the integer $d = \deg L$, and so
\begin{equation*}
\mathcal{B} = \mathcal{B}^{ss} \cup \bigcup_{\substack{d \in \Z \\ d > \frac{1}{2} d_E}} \mathcal{B}_d ,
\end{equation*}
where $\mathcal{B}_d$ is the set of Higgs pairs with Harder-Narasimhan type $d$. For $d>d_E/2$ we define the space $X_d$ to be the union 
\begin{equation}
X_d = \mathcal{B}^{ss} \cup \bigcup_{\substack{\ell \in \Z \\ d \geq \ell > \frac{1}{2} d_E}} \mathcal{B}_\ell \label{eqn:stratification}
\end{equation}
and by convention we set $X_{\lfloor d_E/2\rfloor}= \mathcal{B}^{ss}$.
Then $\{X_d\}_{d=\lfloor d_E/2\rfloor}^\infty$ is the Harder-Narasimhan and YMH-Morse stratification.

This approach for $\mhiggs$ is a special case of a more general method originally outlined by Kirwan, where the topology of a hyperk\"ahler quotient $M \hyperquotient G$ can be studied using a two-step process. First, the cohomology of $\mu_\C^{-1}(0)$ is calculated using the Morse theory of $\left\| \mu_\C \right\|^2$ on $M$ associated to the complex moment map  $\mu_\C= \mu_2+i\mu_3$ , and then the cohomology of $M \hyperquotient G$ can be obtained by studying the K\"ahler quotient of $\mu_\C^{-1}(0)$ by the group $G$ with moment map $\mu_1$. In the case of $M = \totalspace$ we have that $H_\mathcal{G}^*(\totalspace) = H_\mathcal{G}^*(\mathcal{B})$. Therefore, in the Higgs bundle case studied here, it only remains to study the Morse theory of $\YMH = \left\| \mu_1 \right\|^2$ on $\mathcal{B}$ and $\mathcal{B}_0$ respectively.

The formula obtained here for the equivariant cohomology of the minimum has the form
\begin{equation}\label{eqn:generalmorseformula}
P_t^{\mathcal{G}}(B^{ss}) = P_t^{\mathcal{G}}(\mathcal{B}) - \sum_{d=0}^\infty t^{2\mu_d} P_t^{\mathcal{G}}(\mathcal{B}_d) + \sum_{d=1}^{g-1} t^{2\mu_d} P_t^{\mathcal{G}}(\mathcal{B}_{d,\varepsilon}', \mathcal{B}_{d, \varepsilon}'')
\end{equation}
where $\mathcal{B}_d$ denotes the $d^{th}$ stratum of the functional $\YMH$, $\mu_d$ is the rank of a certain bundle over the $d^{th}$ critical set $\eta_d$ (see \eqref{eqn:eta}) representing a subset of the negative eigenspace of the Hessian of $\YMH$ at $\eta_d$, and $P_t^{\mathcal{G}}(\mathcal{B}_{d,\varepsilon}', \mathcal{B}_{d, \varepsilon}'')$ are {\it correction terms} arising from the fact that that the Morse index is not well-defined on the first $g-1$ critical sets. Indeed, as shown in \cite{Wilkin08}, the Morse index at each critical point of $\YMH$ can jump from point to point within the same component of the critical set, and so standard Morse theory cannot be used \emph{a priori}. If the space $\mathcal{B} = \mu_\C^{-1}(0)$ were smooth then the Morse index would be well-defined and the Morse function equivariantly perfect (as is the case for the symplectic reduction considered in \cite{AtiyahBott83} or \cite{Kirwan84}) and the formula for the cohomology of $M \hyperquotient G$ would only consist of the first two terms in \eqref{eqn:generalmorseformula}. 
However, this paper shows that it is possible to construct the Morse theory by hand, using the commutative diagram \eqref{eqn:diagram} in Section \ref{sec:morse}, and computing the cohomology groups of the stratification at each stage.

In order to explain  how to define the index $\mu_d$ in our case we proceed as follows: Regarding $\totalspace$ as the cotangent bundle $T^* \mathcal{A}$, and $\mathcal{B} = \mu_\C^{-1}(0)$ as a subspace of this bundle,  on a critical set of $\YMH$ the solutions of the negative eigenvalue equation of the Hessian of $\YMH = \left\| \mu_1 \right\|^2$ split naturally into two components; one corresponding to the index of the restricted functional $\left\| \left. \mu_1 \right|_\mathcal{A} \right\|^2$, and one along the direction of the cotangent fibers. The dimension of the first component is well-defined over all points of the critical set (this corresponds to $\mu_d$ in the formula above), and the Atiyah-Bott lemma can be applied to the negative normal bundle defined along these directions. The dimension of the second component is not well-defined over all points of the critical set, the methods used here to deal with this show that this leads to extra terms in the Poincar\'e p
 olynomial of $B \mathcal{G}$ corresponding to $P_t^{\mathcal{G}}(\mathcal{B}_{d,\varepsilon}', \mathcal{B}_{d, \varepsilon}'')$. 
More or less this method should work for any hyperk\"ahler quotient of a cotangent bundle.

For the non-fixed determinant case, the long exact sequence obtained at each step of the Morse stratification splits into short exact sequences, thus providing a simple proof of the surjectivity of the hyperk\"ahler Kirwan map. This is done by careful analysis of the correction terms, and it is in a way one of the key observations of this paper (cf.\ Section \ref{sec:surjectivity}). As mentioned above this fails in the  fixed determinant case. 

This paper is organized as follows.
Section \ref{sec:Kuranishi} describes the infinitesimal topology of the stratification arising from the Yang-Mills-Higgs functional.   We define an appropriate linearization of the ``normal bundle'' to the strata and compute its equivariant cohomology.

Section \ref{sec:morse} is  the heart of the paper and contains the details of the Morse theory used to calculate the cohomology of the moduli space. The first result  proves the isomorphism in Proposition  \ref{prop:bott}.  This is the exact analogue of Bott's Lemma
\cite[p.\ 250]{Bott54} in the sense of Bott-Morse theory. The second main result of the section is the commutative diagram \eqref{eqn:diagram} which describes how attaching the strata affects the topology of our space. 
As mentioned before the main difference between Poincare polynomials of hyperk\"ahler quotients from Poincare polynomials of symplectic quotients is the appearance of the rather mysterious correction terms in formula (\ref{eqn:generalmorseformula}). In the course of the proof of Proposition \ref{prop:bott} we show how these 
terms correspond by excision to the fixed points of the $S^1$ action on the moduli space of Higgs bundles. This in our opinion provides an interesting link between our approach and Hitchin's that should be further explored.

Section \ref{sec:surjectivity} contains a detailed analysis of the exact sequence derived from the Morse theory.
We prove Kirwan surjectivity for any degree in the non fixed determinant case (cf.\ Theorem \ref{thm:kirwansurjectivity}). This is achieved by showing that the vertical exact sequence in diagram \eqref{eqn:diagram} splits inducing a splitting on the horizontal sequence. The key to this are results of MacDonald \cite{MacDonald62} on the cohomology of the symmetric product of a curve.
Next, we introduce the fundamental $\Gamma_2 = H^1(M, \mathbb{Z}_2)$ action on the equivariant cohomology which played an important role in the original work of Harder-Narasimhan, Atiyah-Bott and Hitchin (cf.\ \cite{AtiyahBott83, Hitchin87}). The action  splits the exact sequences in diagram \eqref{eqn:diagram} into $\Gamma_2$-invariant and noninvariant parts, and the main result is Theorem \ref{thm:invariantkirwansurjectivity}, which demonstrates Kirwan surjectivity holds on $\Gamma_2$-invariant part of the cohomology.

Finally, Section \ref{sec:computations} contains the computations of the equivariant Poincar\'e polynomials of $\mathcal{B}^{ss}$ and $\mathcal{B}_0^{ss}$ stated above. \\

{\bf Acknowledgments.} We are thankful to Megumi Harada, Nan-Kuo Ho and Melissa Liu for pointing out an error in a previous version of the paper.

\section{Local structure of the space of Higgs bundles} \label{sec:Kuranishi}

In this section we explain the Kuranishi model for Higgs bundles (cf.\ \cite{ahs78} and \cite[Ch.\ VII]{Kobayashi87}) and   derive the basic results needed for the Morse theory of Section \ref{sec:morse}. For simplicity, we treat the case of non-fixed determinant, and the results for fixed determinant are identical \emph{mutatis mutandi}.

\subsection{The deformation complex}
We begin with the deformation theory.

Infinitesimal deformations of $(A, \Phi) \in \mathcal{B}$ modulo equivalence are described by the following elliptic complex, which we denote by $\mathcal{C}_{(A, \Phi)}$.

\begin{align} 
\begin{split}
\xymatrix{
\mathcal{C}_{(A, \Phi)}^0 \ar[r]^(0.45){D_1} &\mathcal{C}_{(A, \Phi)}^1\ar[r]^(0.45){D_2} & \mathcal{C}_{(A, \Phi)}^2
} \hskip1.25in
\label{eqn:deformation-complex}\\
\xymatrix{
\Omega^0(\End(E)) \ar[r]^(0.35){D_1} & \Omega^{0,1}(\End(E)) \oplus \Omega^{1,0}(\End(E)) \ar[r]^(0.65){D_2} & \Omega^2(\End(E)) 
}
\end{split}
\end{align}
where 
$$
D_1(u)  = (d_A'' u, [\Phi, u])
\ ,\
D_2(a, \varphi)  = d_A'' \varphi + [a, \Phi]
$$
Here, $D_1$ is the linearization of the action of the complex gauge group on $\mathcal{B}$, and $D_2$ is the linearization of the condition $d_A''\Phi=0$.
Note that $D_2 D_1 = [d_A'' \Phi, u] = 0$ if $(A, \Phi) \in \mathcal{B}$. 

The hermitian metric gives adjoint operators $D_1^\ast$, $D_2^\ast$, and the spaces of \emph{harmonic forms} are given by 
\begin{align*}
{\mathcal H}^0(\mathcal{C}_{(A, \Phi)})&=\ker D_1 \\
{\mathcal H}^1(\mathcal{C}_{(A, \Phi)})&=\ker D_1^\ast\cap \ker D_2 \\
{\mathcal H}^2(\mathcal{C}_{(A, \Phi)})&=\ker D_2^\ast 
\end{align*}
with harmonic projections $\Pi_i: \mathcal{C}_{(A, \Phi)}^i\to {\mathcal H}^i(\mathcal{C}_{(A, \Phi)})$.

We will be interested in the deformation complex along higher critical sets of the Yang-Mills-Higgs functional .  These are given by split Higgs bundles $(A,\Phi)=(A_1\oplus A_2, \Phi_1\oplus \Phi_2)$ corresponding to a smooth splitting 
 $E = L_1 \oplus L_2$ of $E$ with $\deg L_1=d>\deg L_2=d_E-d$.  The set of all such critical points is denoted by $\eta_d\subset\B$.   We will often use the notation $L=L_1\otimes L_2^\ast$, and $\Phi_\flat=\tfrac{1}{2}(\Phi_1-\Phi_2)$, and denote the  components of $\End(E)\simeq L_i\otimes L_j^\ast$ in the complex by $u_{ij}$, $a_{ij}$, $\varphi_{ij}$, $u_\flat=\tfrac{1}{2}(u_{11}-u_{22})$,  etc.
Define $\End(E)^{UT}$ to be the subbundle of $\End(E)$ consisting of endomorphisms that preserve $L_1$, and  $\End(E)^{SUT} \subset \End(E)^{UT}$ to be the subbundle of endomorphisms whose component in the subbundle $\End(L_1) \oplus \End(L_2)$ is zero. We say that 
$$(a, \varphi) \in \Omega^{0,1}(\End(E)^{UT}) \oplus \Omega^{1,0}(\End(E)^{UT})$$ is \emph{upper-triangular}, and $$(a, \varphi) \in \Omega^{0,1}(\End(E)^{SUT}) \oplus \Omega^{1,0}(\End(E)^{SUT})$$ is \emph{strictly upper-triangular}. Similarly, define the \emph{lower-triangular}, \emph{strictly lower-triangular},  \emph{diagonal} and \emph{off-diagonal}   endomorphisms, with the obvious notation.  
Since $\Phi$ is diagonal, harmonic projection  preserves components.   For example, ${\mathcal H}^1({\mathcal C}_{(A,\Phi)})$ consists of all $(a,\varphi)$ satisfying
\begin{eqnarray}
d''\varphi_{ii}=0 & (d'')^\ast a_{ii}=0 \label{eqn:harmonic-diag} \\
d_A''\varphi_{12}+2\Phi_\flat a_{12}=0 & (d_A'')^\ast a_{12}+2\bar\ast(\Phi_{\flat}\bar\ast\varphi_{12})=0 \label{eqn:harmonic-12} \\
d_A''\varphi_{21}-2\Phi_{\flat}a_{21}=0 &
(d_A'')^\ast a_{21}-2\bar\ast(\Phi_{\flat}\bar\ast\varphi_{21})=0 \label{eqn:harmonic-21}\end{eqnarray}
where $\bar\ast$ is defined as in \cite[eq.\ (2.8)]{Kobayashi87}. 

The following construction will be important for the computations in this paper.
\begin{defi}\label{def:negative-bundles}
Let $\tilde{q} : \mathcal{T} \rightarrow \eta_{d}$ be the trivial bundle over
 $\eta_{d}$ with fiber 
 $$\Omega^{0,1}(\End(E)) \oplus \Omega^{1,0}(\End(E))$$
  and define $\nu_d^- \subset \mathcal{T}$ to be the subspace with projection map $\tilde{q} : \nu_d^- 
  \rightarrow \eta_{d}$, where the fiber over $(A,\Phi) \in \eta_{d}$ is ${\mathcal H}^1({\mathcal C}_{(A,\Phi)}^{SLT})$.
   Note that in general the 
  dimension of the fiber may depend on the Higgs structure.

We also define the subsets
\begin{align*}
\nu_d' &= \nu_d^- \setminus \eta_{d} \\
\nu_d'' &= \left\{ ((A,\Phi), (a, \varphi)) \in \nu_d^- \, : \, {\mathcal H}(a_{21}) \neq 0 \right\} 
\end{align*}
 where $\mathcal H$ denotes the $d_A''$-harmonic projection.
\end{defi}

\subsection{Equivariant cohomology of the normal spaces}
Note that there is a natural action of $\G$ on the spaces introduced in Definition \ref{def:negative-bundles}.  In this section we compute the $\G$-equivariant cohomology associated to the triple $\nu_d^-$, $\nu_d'$, and $\nu_d''$.
We first  make the following
 \begin{defi} \label{def:psi}
 Let $(A,\Phi)\in \eta_d$, $E=L_1\oplus L_2$, and  $L=L_1\otimes L_2^\ast$.
  Let $(a,\varphi)\in \Omega^{0,1}(\End(E))\oplus \Omega^{1,0}(\End(E))$. Since $\deg L>0$,  there is a unique $f_{21}\in \Omega^0(L^\ast)$ such that
  $a_{21}={\mathcal H}(a_{21})+d_A''f_{21}$. 
 Define
 \begin{equation} \label{eqn:psi}
 \Psi:\Omega^{0,1}(\End(E))\oplus \Omega^{1,0}(\End(E))\to H^{1,0}(L): (a,\varphi)\mapsto{\mathcal H}(\varphi_{21}+2f_{21}\Phi_{\flat})
 \end{equation}
  Set $\psi_{21}=\varphi_{21}+2f_{21}\Phi_{\flat}$, and let $F_{21}$ be the unique section in $(\ker (d_A'')^\ast)^\perp\subset \Omega^{1,1}(L^\ast)$ such that $\psi_{21}=\Psi(a,\varphi)+(d_A'')^\ast F_{21}$.
 \end{defi}
Let
\begin{equation} 
T_d=\left\{ (a,\varphi)\in \nu_d^- : {\mathcal H}(a_{21})=0\ ,\ \Psi(a,\varphi)\neq 0\right\}\label{eqn:td}
\end{equation}
and set $\mu_d=g-1+2d-d_E$.
We will prove the following
\begin{thm} \label{thm:nucohomology}
There are isomorphisms
\begin{align}
H_\G^\ast(\nu_d^-,\nu_d'')&\simeq H^{\ast-2\mu_d}_\G(\eta_{d})
\label{eqn:nucohomology-doubleprime} \\
H_\G^\ast(\nu_d',\nu_d'')&\simeq H^{\ast-2\mu_d}_\G(T_d)
\label{eqn:nucohomology-prime-doubleprime}
\end{align}
\end{thm}

With the notation above,
  eq.\ \eqref{eqn:harmonic-21} becomes
 \begin{align}
 \square_A F_{21}+2{\mathcal H}(a_{21})\Phi_{\flat}&= 0 \label{eqn:F21} \\
    \square_A f_{21}+4\Vert\Phi_{\flat}\Vert^2f_{21}&=(d_A'')^\ast (2\bar\ast(\Phi_{\flat}\bar\ast F_{21})
)+ 2\bar\ast(\Phi_{\flat}\bar\ast\Psi(a,\varphi))
 \label{eqn:f21}
 \end{align}
  Given
 $({\mathcal H}(a_{21}), \Psi(a,\varphi))\in H^{0,1}(L^\ast)\oplus H^{1,0}(L^\ast)$,
  satisfying  ${\mathcal H}({\mathcal H}(a_{21})\Phi_{\flat})=0$, \eqref{eqn:F21} uniquely determines $F_{21}$. 
  Then \eqref{eqn:f21} uniquely determines $f_{21}$.  Note that since $\deg L>0$, $\square_A$, and therefore $\square_A +\Vert\Phi_{\flat}\Vert^2$, has no kernel. 
   We then reconstruct $(a,\varphi)\in {\mathcal H}^1({\mathcal C}_{(A,\Phi)}^{SLT})$ by setting
  $a_{21}={\mathcal H}(a_{21})+d_A''f_{21}$, and $\varphi_{21}=\Psi(a,\varphi)+(d_A'')^\ast F_{21}-2f_{21}\Phi_{\flat}$.
 Thus, we have shown
 \begin{align}
 \begin{split} \label{eqn:nu-fiber}
 &\nu_d^-\cap \,\tilde q^{-1}(A,\Phi)\\
 &\simeq
 \left\{ ({\mathcal H}(a_{21}), \Psi(a,\varphi))\in H^{0,1}(L^\ast)\oplus H^{1,0}(L^\ast) :
 {\mathcal H}({\mathcal H}(a_{21})\Phi_{\flat})=0\right\}
 \end{split}
 \end{align}
Next, let $\eta_{d,0}\subset \eta_{d}$ be the subset of critical points where $\Phi=0$.  Notice that $\eta_{d,0}\hookrightarrow \eta_{d}$ is a $\G$-equivariant deformation retraction under scaling $(A,\Phi)\mapsto (A,t\Phi)$, for $0\leq t\leq 1$.   Let
$$
\nu_{d,0}^-=\nu_d^-\cap\, \tilde q^{-1}(\eta_{d,0})\ , \
\nu_{d,0}'=\nu_d'\cap\, \tilde q^{-1}(\eta_{d,0})\ , \
\nu_{d,0}''=\nu_d''\cap\, \tilde q^{-1}(\eta_{d,0})
$$
 We have the following
\begin{lem} \label{lem:nu-retraction}
There is a $\G$-equivariant retraction $\nu_{d,0}^-\hookrightarrow\nu_d^-$ that preserves the subspaces
$\nu_d'$ and $\nu_d''$.
\end{lem}

 \begin{proof}
 Given $(A,\Phi)$ and $(a_{21},\varphi_{21})\in \Omega^{0,1}(L^\ast)\oplus \Omega^{1,0}(L^\ast)$, let $(f_{21}(\Phi), F_{21}(\Phi))$ be the unique solutions to \eqref{eqn:F21} and \eqref{eqn:f21}.  Notice that $(f_{21}(0), F_{21}(0))=(0,0)$.
  Then an explicit retraction may be defined as follows
 \begin{align*}
 \rho : [0,1]\times \nu_d^-&\longrightarrow \nu_d^-  \\
\rho (t, (A,\Phi), (a,\varphi)) &= \bigl( (A,t\Phi), {\mathcal H}(a_{21})+d_A''f_{21}(t\Phi), \\
&\qquad\qquad \Psi(a,\varphi)+(d_A'')^\ast F_{21}(t\Phi)-2tf_{21}(t\Phi)\Phi_{\flat} \bigr)
 \end{align*}
 It is easily verified that $\rho$ satisfies the properties stated in the lemma.
 \end{proof}
 
 \begin{proof}[Proof of Theorem \ref{thm:nucohomology}]
First, note that by Riemann-Roch, $\dim H^{0,1}(L^\ast)=\mu_d$. By Lemma \ref{lem:nu-retraction},
there are $\G$-equivariant homotopy equivalences $(\nu_{d,0}^-,\nu_{d,0}'')\simeq (\nu_{d}^-,\nu_{d}'')$, and $(\nu_{d,0}',\nu_{d,0}'')\simeq (\nu_{d}',\nu_{d}'')$.
Also, since $\eta_{d,0}\hookrightarrow \eta_{d}$ is a $\G$-equivariant deformation retraction, $H^{\ast}_\G(\eta_{d})\simeq H^{\ast}_\G(\eta_{d,0})$.   By \eqref{eqn:nu-fiber}, a similar statement holds for
\begin{equation}\label{eqn:td0}
T_{d,0}=T_d\cap \tilde q^{-1}(\eta_{d,0})
\end{equation}
Hence, it suffices to prove 
\begin{align}
H_\G^\ast(\nu_{d,0}^-,\nu_{d,0}'')&\simeq H^{\ast-2\mu_d}_\G(\eta_{d,0})
\label{eqn:nu0-doubleprime} \\
H_\G^\ast(\nu_{d,0}',\nu_{d,0}'')&\simeq H^{\ast-2\mu_d}_\G(T_{d,0})
\label{eqn:nu0-prime-doubleprime}
\end{align}
From \eqref{eqn:nu-fiber} we have,
  \begin{equation*} \label{eqn:nu0-fiber}
 \nu_{d,0}^-\,\cap\tilde q^{-1}(A,0) \simeq
  H^{0,1}(L^\ast)\oplus H^{1,0}(L^\ast) 
 \end{equation*}
 Then \eqref{eqn:nu0-doubleprime} follows from 
this and the Thom isomorphism theorem.  Next, let
$$
Y_d=\left\{ (a,\varphi)\in \nu_{d,0}'' : \Psi(a,\varphi)=0\right\}
$$
Clearly, $Y_d$ is closed in $\nu_{d,0}''$, and one observes that it is also closed in  $\nu_{d,0}'$.  

 Hence, by excision and the Thom isomorphism applied to the projection to $H^{1,0}(L^\ast)$,
$$
H_\G^\ast(\nu_{d,0}',\nu_{d,0}'')\simeq H_\G^\ast(\nu_{d,0}'\setminus Y_d,\nu_{d,0}''\setminus Y_d)\simeq H_\G^{\ast-2\mu_d}(T_{d,0})
$$
This proves \eqref{eqn:nu0-prime-doubleprime}.
\end{proof}
 
There is an important connection between the topology of the space $T_{d,0}$
 and the fixed points of the $S^1$-action on the
 moduli space of semistable Higgs bundles, and this will be used below.
 Recall from 
\cite[Sec.\ 7]{Hitchin87}
that the non-minimal critical point set of the function $\Vert \Phi\Vert^2$ 
on ${\mathcal M}^{Higgs}(2,d_E)$
 has components
$c_d$ corresponding to equivalence classes of (stable) 
Higgs pairs $(A, \Phi)$, where $A=A_1\oplus A_2$ is a split connection on  
$E = L_1 \oplus L_2$ with $\deg L_1=d>\deg L_2=d_E-d$ and 
$\Phi \neq 0$ is strictly lower triangular with respect to the splitting. 
On the other hand, it follows from  \eqref{eqn:harmonic-21}  and \eqref{eqn:td0}
that 
$$T_{d,0}=\left\{ 
((A, \Phi=0), (\alpha_{21}=0, \varphi_{21})) :
A=A_1\oplus A_2\ , \
d_A''(\varphi_{21})=0\right\}
$$
Taking into account gauge equivalence,  we therefore obtain the following
\begin{lem} \label{music}
Let $c_d$ be as above. For the non-fixed 
determinant case,
 $$H^\ast_\G(T_{d,0})=H^\ast(c_d) \otimes H^*(BU(1))$$
and 
in the fixed determinant case, $H^\ast_\G(T_{d,0})=H^\ast(c_d)$.
\end{lem}

\section{Morse Theory on the space of Higgs bundles}\label{sec:morse}

The purpose of this section is to derive the theoretical results underpinning the calculations in Section \ref{sec:computations}. This is done in a natural way, using the functional $\YMH$ as a Morse function on the singular space $\mathcal{B}$. As a consequence, we obtain a criterion for hyperk\"ahler Kirwan surjectivity in Corollary \ref{cor:surjective}, which we show is satisfied for the non-fixed determinant case in Section \ref{sec:surjectivity}. The key steps in this process are (a) the proof of the isomorphism \eqref{eqn:bott}, which relates the topology of a neighborhood of the stratum to the topology of the negative eigenspace of the Hessian on the critical set (a generalization of Bott's isomorphism  \cite[p.\ 250]{Bott54} to the singular space of Higgs bundles), and (b) the commutative diagram \eqref{eqn:diagram}, which provides a way to measure the imperfections of the Morse function $\YMH$ caused by the singularities in the space $\mathcal{B}$. 

The methods of this section are also valid for the rank $2$ degree $1$ case, and in Section \ref{sec:computations} they are used provide new computations of the results of \cite{Hitchin87} (fixed determinant case) and \cite{HauselThaddeus04} (non-fixed determinant case).


\subsection{Relationship to Morse-Bott theory}
  Recall the spaces  $\nu_d^-$,  $\nu_d'$ and $\nu_d''$ from Definition \ref{def:negative-bundles}. 
 This section is devoted to the proof of the Bott isomorphism
 \begin{prop} \label{prop:bott} For $d>d_E/2$, there is an isomorphism
\begin{equation} \label{eqn:bott}
 H^\ast_\G(X_d,X_{d-1})\simeq H^\ast_{\G}(\nu_d^-, \nu_d')
 \end{equation}
 \end{prop}
  Let ${\mathcal A}_d$ denote
 the stable manifold in $\mathcal A$ of the critical set $\eta_{d,0}$ of the Yang-Mills functional (cf.\  \cite{AtiyahBott83, Daskal92}). We also
 define 
 $$
 X_d^{\mathcal A}={\mathcal A}^{ss}\cup \bigcup_{d_E/2<\ell\leq d} {\mathcal A}_\ell
 $$
 Let $X_d''=X_d\setminus\pr^{-1}({\mathcal A}_d)$.   By applying the five lemma to the exact sequences for the triples $(X_d, X_{d-1}, X_d'')$ and $(\nu_d^-,  \nu_d', \nu_d'')$, it suffices to prove the two isomorphisms
  \begin{align}
   H^\ast_\G(X_d, X_d^{\prime\prime})&\simeq H^\ast_{\G}(\nu_d^-, \nu_d'') \label{eqn:doubleprime} \\
   H^\ast_\G(X_{d-1},X_{d}'')&\simeq H^\ast_{\G}(\nu_d', \nu_d'') \ \label{eqn:prime-doubleprime} .
  \end{align}

We begin with the first equality.
 \begin{proof}[Proof of \eqref{eqn:doubleprime}]
 By  \eqref{eqn:nucohomology-doubleprime}, the result of Atiyah-Bott \cite{AtiyahBott83}, and
 the fact that the projection $\eta_{d}\to \eta_{d,0}$ has contractible fibers, it suffices to show
 $$ H^\ast_\G(X_d, X_d^{\prime\prime})\simeq H^\ast_{\G}(X_d^{\mathcal A}, X_{d-1}^{\mathcal A})$$
 Also, note that for $\ell> d/2$, $\pr(\B_\ell)={\mathcal A}_\ell$.  Indeed, the inclusion $\supset$ comes from taking $\Phi=0$, and the inclusion $\subset$ follows from the fact that for any extension of line bundles
 $$
 0\longrightarrow L_1\longrightarrow E\longrightarrow L_2\longrightarrow 0
 $$
 with $\deg L_1>\deg L_2$, $0\subset L_1\subset E$ is precisely the Harder-Narasimhan filtration of $E$.  With this understood, let
$
{\mathcal K}_d=\pr(\B^{ss})\cap\left(\cup_{\ell >d}{\mathcal A}_\ell\right)
$.  Then 
we claim that ${\mathcal K}_d$, which is manifestly contained in $\pr(X_d)$,  is in fact closed in $\pr(X_d)$.  To see this, let $A_j\in {\mathcal K}_d$, $A_j\to A\in \pr(X_d)$.  By definition, $A=\pr(A,\Phi)$ with either $(A,\Phi)\in \B^{ss}$, or $(A,\Phi)\in \B_\ell$, $\ell\leq d$.  Notice that by semicontinuity, $A\in \cup_{\ell >d}{\mathcal A}_\ell$.  Hence, the second possibility does not occur.
It must therefore be the case that $A\in \pr(\B^{ss})$, and hence $A\in {\mathcal K}_d$ also.  Now, since ${\mathcal K}_d\cap{\mathcal A}_d=\emptyset$ by definition, it follows that
$$
{\mathcal K}_d\subset \pr(X_d)\setminus {\mathcal A}_d=\pr(X_d\setminus \pr^{-1}({\mathcal A}_d))=\pr(X_d'')
$$
 Since the fibers of the map $\pr:X_d\to \pr(X_d)$ are $\G$-equivariantly  contractible via scaling of the Higgs field,  it follows from excision that
 $$
 H^\ast_\G(X_d, X_d^{\prime\prime})\simeq H^\ast_\G(\pr(X_d), \pr(X_d^{\prime\prime}))\simeq
 H^\ast_\G(\pr(X_d)\setminus{\mathcal K}_d, \pr(X_d^{\prime\prime})\setminus{\mathcal K}_d)
 $$
 However,
 \begin{align*}
 \pr(X_d)&=\pr(\B^{ss})\cup \left(\cup_{d_E/2<\ell\leq d}  \pr(\B_\ell)\right) \\
 &={\mathcal A}^{ss}\cup\left(\cup_{ d_E/2<\ell}{\mathcal A}_{\ell}\cap \pr(\B^{ss}) \right)\cup \left(\cup_{d_E/2<\ell\leq d} {\mathcal A}_{\ell}\right) \\
 &={\mathcal A}^{ss}\cup {\mathcal K}_d\cup\left( \cup_{d_E/2<\ell\leq d} {\mathcal A}_{\ell} \right)
 \end{align*}
 Hence, since the union is disjoint, $\pr(X_d)\setminus {\mathcal K}_d=X_d^{\mathcal A}$.  Furthermore, 
 $$
 \pr(X_d'')\setminus {\mathcal K}_d=\pr(X_d)\setminus {\mathcal K}_d\cup {\mathcal A}_d=X_d^{\mathcal A}\setminus {\mathcal A}_d=X_{d-1}^{\mathcal A}
 $$
 This completes the proof.
 \end{proof}

\begin{proof}[Proof of \eqref{eqn:prime-doubleprime}]
By the isomorphism \eqref{eqn:nucohomology-prime-doubleprime} (see also Lemma \ref{lem:nu-retraction}), it suffices to prove $H^\ast_{\G}(X_{d-1},X_d'')\simeq H^\ast_{\G}(T_{d,0})$.  From the proof of \eqref{eqn:doubleprime} we have
\begin{align*}
X_d''&=\left\{\B^{ss}\cup(\cup_{d_E/2<\ell\leq d} \B_\ell)\right\}\setminus \pr^{-1}({\mathcal A}_d) \\
&=\left\{\B^{ss}\setminus \pr^{-1}({\mathcal A}_d)\right\}\cup(\cup_{d_E/2<\ell\leq d-1} \B_\ell)
\end{align*}
whereas 
$$
X_{d-1}=\B^{ss}\cup(\cup_{d_E/2<\ell\leq d-1} \B_\ell)
$$
Since $\cup_{d_E/2<\ell\leq d-1} \B_\ell\subset X_d''$ is closed in $X_{d-1}$, it follows from excision that
$$
H^\ast_\G(X_{d-1}, X_d'')\simeq H^\ast_\G(\B^{ss},\B^{ss}\setminus \pr^{-1}({\mathcal A}_d))
$$
By the main result of \cite{Wilkin08}, the YMH-flow gives  a $\G$-equivariant deformation retract to $\B_{min}$.  Hence,
$$
H^\ast_\G(X_{d-1}, X_d'')\simeq H^\ast_\G(\B_{min},\B_{min}\setminus \pr^{-1}({\mathcal A}_d) ).
$$
Next, notice that the singularities of $\B_{min}$ correspond to strictly semistable points and therefore there exists a neighborhood 
$\mathcal N_d$ of  $\pr^{-1}({\mathcal A}_d) \cap \B_{min}$ in $\B_{min}$ consisting entirely of smooth points.  Furthermore, $\G$ acts on 
$\mathcal N_d$ with constant central transformations as stabilizers. Therefore, by again applying excision and passing to the quotient we obtain 
\begin{eqnarray*}
H^\ast_\G(X_{d-1}, X_d'') &\simeq& H^\ast_\G(\mathcal N_d,  \mathcal N_d \setminus \pr^{-1}({\mathcal A}_d)) \\
&\simeq& H^\ast (\mathcal N_d \slash \G ,  (\mathcal N_d \setminus \pr^{-1}({\mathcal A}_d)) \slash \G)\otimes H^{\ast}(BU(1)).
\end{eqnarray*}
Now according to Frankel and Hitchin (cf.\ \cite[Sect.\ 7]{Hitchin87}) the latter equality localizes the computation to the $d$-th component $c_d$ of the fixed point set for the $S^1$-action on $\B_{min}/\G$.
Hence,
$$H^\ast_\G (X_{d-1}, X_d'') \simeq  H^{\ast-2\mu_d} ( c_d) \otimes H^{\ast} (BU(1)).$$
The result follows by combining the above isomorphism with Theorem \ref{thm:nucohomology} and Lemma \ref{music}.
\end{proof}

\subsection{A framework for cohomology computations}
From Proposition \ref{prop:bott}, the computation of $H_\mathcal{G}^*(\nu_d^-, \nu_d')$ in Theorem \ref{thm:nucohomology}  leads to a computation of the equivariant cohomology of  the space of rank $2$ Higgs bundles, using the commutative diagram \eqref{eqn:diagram}. Recall the 
decomposition \eqref{eqn:stratification}.

The inclusion $X_{d-1} \hookrightarrow X_d$ induces a long exact sequence in equivariant cohomology
\begin{equation}\label{eqn:LES-for-Xd}
\cdots \rightarrow H_\mathcal{G}^*(X_d, X_{d-1}) \rightarrow H_\mathcal{G}^* (X_d) \rightarrow H_\mathcal{G}^*(X_{d-1}) \rightarrow \cdots ,
\end{equation}
and the method of this section is to relate the cohomology groups $H_\mathcal{G}^* (X_d)$ and $H_\mathcal{G}^*(X_{d-1})$ by $H_\mathcal{G}^*(X_d, X_{d-1})$ and the maps in the corresponding long exact sequence for $(\nu_d^-, \nu_d', \nu_d'')$.

Let $J_d(M)$ denote the Jacobian of degree $d$ line bundles over the Riemann surface $M$, let $S^n M$ denote the $n^{th}$ symmetric product of $M$, and let $\widetilde{S}^n M$ denote the $2^{2g}$ cover of $S^n M$ described in   \cite[eq.\ (7.10)]{Hitchin87}. The critical sets correspond to $\Phi$-invariant holomorphic splittings $E = L_1 \oplus L_2$, therefore after dividing by the unitary gauge group $\mathcal{G}$ the critical sets of $\YMH$ are
\begin{equation} \label{eqn:eta}
\eta_d=\begin{cases}
 T^* J_d(M) \times T^*J_{d_E-d}(M) & \text{non-fixed determinant case;}\\
 T^* J_d(M) &\text{fixed determinant case.}
\end{cases}
\end{equation}

By combining this with Lemma \ref{music} and the computation in \cite{Hitchin87} we obtain
\begin{lem}\label{lem:critsetcohomology}
In the non-fixed determinant case
\begin{align}
H_\mathcal{G}^*(\eta_d) & \cong H^*(J_d(M) \times J_n(M)) \otimes H^*(BU(1))^{\otimes 2} \label{eqn:nonfixedcohomologycritical} \\
H_\mathcal{G}^*(T_d) & \cong H^*(J_d(M)) \otimes H^*(S^n M) \otimes H^*(BU(1)) \label{eqn:nonfixednucohomology-prime-doubleprime} .
\end{align}
In the fixed determinant case
\begin{align}
H_\mathcal{G}^*(\eta_d) & \cong H^*(J_d(M)) \otimes H^*(BU(1)) \label{eqn:fixedcohomologycritical} \\
H_\mathcal{G}^*(T_d) & \cong H^*(\widetilde{S}^n M) \label{eqn:fixednucohomology-prime-doubleprime} .
\end{align}
\end{lem}

The spaces $(\nu_d^-, \nu_d', \nu_d'')$ form a triple, and the isomorphism $H_\mathcal{G}^*(X_d, X_{d-1}) \cong H_\mathcal{G}^* (\nu_d^-, \nu_d')$ from \eqref{eqn:bott} implies the long exact sequence (abbrev.\ \emph{LES}) of this triple is related to the LES \eqref{eqn:LES-for-Xd} in the following commutative diagram.
\begin{equation}\label{eqn:diagram}
\xymatrix{
& \vdots \ar[d]^{\delta^{k-1}} \\
\cdots \ar[r] & H_{\mathcal{G}}^k(X_d, X_{d-1}) \ar[d]^{\cong} \ar[r]^(0.6){\alpha^k} & H_{\mathcal{G}}^k(X_d) \ar[r]^{\beta^k} \ar[d] & H_{\mathcal{G}}^k(X_{d-1}) \ar[d] \ar[r]^(0.6){\gamma^k} & \cdots \\
\cdots \ar[r] & H_\mathcal{G}^k(\nu_d^-, \nu_d') \ar[d]^{\zeta^k} \ar[r]^(0.6){\alpha_\varepsilon^k} & H_\mathcal{G}^k(\nu_d^-) \ar[r]^{\beta_\varepsilon^k} \ar[d]^{\omega^k} & H_\mathcal{G}^k(\nu_d') \ar[r]^(0.6){\gamma_\varepsilon^k} & \cdots \\
& H_\mathcal{G}^k(\nu_d^-, \nu_d'') \ar[r]^(0.6){\smallsmile e} \ar[d]^{\lambda^k} \ar[ur]^{\xi^k} & H_\mathcal{G}^k(\eta_d) \\
& H_\mathcal{G}^k(\nu_d', \nu_d'') \ar[d]^{\delta^k} \\
& \vdots \\
}
\end{equation}
where the two horizontal exact sequences are the LES of the pairs $(X_d, X_{d-1})$ and $(\nu_d^-, \nu_d')$ respectively. The vertical exact sequence in the diagram is the LES of the triple $(\nu_d^-, \nu_d', \nu_d'')$. The diagonal map $\xi^k$ is from the LES of the pair $(\nu_d^-, \nu_d'')$. Applying the Atiyah-Bott lemma (\cite[Prop.\ 13.4]{AtiyahBott83}) gives us the following lemma.
\begin{lem}\label{lem:Eulerclass}
The map $\smallsmile e : H_{\mathcal{G}}^k(\nu_d^-, \nu_d'') \rightarrow H_{\mathcal{G}}^k(\eta_d)$ is injective and therefore the map $\xi^k$ is injective, since $\omega^k \circ \xi^k = \smallsmile e$.
\end{lem} 
From the horizontal LES of (\ref{eqn:diagram})
\begin{equation*}
\frac{H_{\mathcal{G}}^k(X_{d-1})}{\im \beta^k} \cong \frac{H_{\mathcal{G}}^k(X_{d-1})}{\ker \gamma^k} \cong \im \gamma^k \cong \ker \alpha^{k+1}
\end{equation*}
and also
\begin{equation*}
\im \beta^k \cong \frac{H_{\mathcal{G}}^k(X_d)}{\ker \beta^k} \cong \frac{H_{\mathcal{G}}^k(X_d)}{\im \alpha^k}
\end{equation*}
Therefore
\begin{align*}
\dim \, \ker \alpha^{k+1} & = \dim H_{\mathcal{G}}^k(X_{d-1}) - \dim \, \im \beta^k \\
 & = \dim H_{\mathcal{G}}^k(X_{d-1}) - \dim H_{\mathcal{G}}^k(X_d) + \dim \, \im \alpha^k
\end{align*}

\begin{lem}\label{lem:subset}

$\ker \alpha^k \subseteq \ker \zeta^k$.

\end{lem}

\begin{proof}

Lemma \ref{lem:Eulerclass} implies $\xi^k$ is injective, and since $\alpha_\varepsilon^k = \xi^k \circ \zeta^k$, then $\ker \alpha_\varepsilon^k = \ker \zeta^k$. Using the isomorphism \eqref{eqn:bott} to identify the spaces $H_\mathcal{G}^*(X_d, X_{d-1}) \cong H_\mathcal{G}^*(\nu_d^-, \nu_d')$, we see that $\ker \alpha^k \subseteq \ker \alpha_\varepsilon^k$, which completes the proof.
\end{proof}

\begin{cor}\label{cor:surjective}
If $\lambda^k$ is surjective for all $k$, then $\beta^k$ is surjective for all $k$.
\end{cor}

\begin{proof}
If $\lambda^k$ is surjective for all $k$, then $\zeta^k$ is injective for all $k$, and so Lemma \ref{lem:subset} implies $\alpha^k$ is injective for all $k$. Therefore, $\beta^k$ is surjective for all $k$.
\end{proof}

In particular, we see that if for each stratum $X_d$, we can show that $\lambda^k$ is surjective for all $k$, then the inclusion $\mathcal{B}^{ss} \hookrightarrow \mathcal{B}$ induces a surjective map $\kappa_H : H_\mathcal{G}^*(\mathcal{B}) \rightarrow H_\mathcal{G}^*(\mathcal{B}^{ss})$. The next section shows that this is indeed the case for non-fixed determinant Higgs bundles.

\section{Hyperk\"ahler Kirwan surjectivity}\label{sec:hksurjectivity}

We now apply the results of  Section \ref{sec:morse} to the question of Kirwan surjectivity for Higgs bundles.  We establish surjectivity in the case of the non-fixed determinant moduli space.  In the fixed determinant case surjectivity fails, this will be explained in more detail in Section \ref{sec:fixedsurjectivity}, where we introduce an action of $\Gamma_2 = H^1(M, \mathbb{Z}_2)$ and prove surjectivity onto the $\Gamma_2$-invariant equivariant cohomology.

\subsection{The non-fixed determinant case}  \label{sec:surjectivity}

For simplicity of notation, throughout this section let $n = 2g - 2 +d_E - 2d$ where $d_E = \deg (E)$ and $d$ is the index of the stratum $\mathcal{B}_d$ as defined in Section \ref{sec:morse}. In this section we prove

\begin{thm}\label{thm:kirwansurjectivity}
The spaces $\mhiggs(2, 1)$ and $\mhiggs(2, 0)$ are hyperk\"ahler quotients $T^* \mathcal{A} \hyperquotient \mathcal{G}$ for which the hyperk\"ahler Kirwan map $$\kappa_{H} : H_\mathcal{G}^*(T^* \mathcal{A}) \rightarrow H_\mathcal{G}^*(\mathcal{B}^{ss})$$ is surjective.
\end{thm}

As mentioned in the Introduction, for the space $\mhiggs(2, 1)$ a special case of Theorem \ref{thm:kirwansurjectivity} has already been proven by Hausel and Thaddeus in \cite{HauselThaddeus04}. However, because of singularities  their methods do not apply to the space $\mhiggs(2,0)$. 

The calculations of Hitchin in \cite{Hitchin87} for $\mhiggs_0(2,1)$, and those of Section \ref{sec:computations} in this paper for $\mhiggs_0(2,0)$, show that the hyperk\"ahler Kirwan map cannot be surjective for the fixed determinant case. The results of this section also provide a basis for the proof of Theorem \ref{thm:invariantkirwansurjectivity} below, where we show that the hyperk\"ahler Kirwan map is surjective onto the $\Gamma_2$-invariant part of the cohomology.  This is the best possible result for the fixed determinant case. 

The proof of Theorem \ref{thm:kirwansurjectivity} reduces to showing that the LES \eqref{eqn:LES-for-Xd} splits, and hence the map $\beta^* : H_\mathcal{G}^*(X_d) \rightarrow H_\mathcal{G}^*(X_{d-1})$ is surjective for each positive integer $d$. 
Lemma \ref{lem:subset} shows that this is the case iff the vertical LES of diagram \eqref{eqn:diagram} splits. By Corollary \ref{cor:surjective}, together with the description of the cohomology groups in Theorem \ref{thm:nucohomology},  the proof of Theorem \ref{thm:kirwansurjectivity} reduces to showing that the map $\lambda^* : H_\mathcal{G}^{*-2\mu_d}(\eta_d) \rightarrow H_\mathcal{G}^{*-2\mu_d}(T_d)$ is surjective. In the non-fixed determinant case, the following lemma provides a simpler description of the map $\lambda^*$.

\begin{lem}\label{cla:restrictedsurjective}
 The map $\lambda^*$ restricts to a map 
\begin{equation*}
\lambda_r^* : H^{*-2\mu_d} (J_d(M)) \otimes H^*(BU(1)) \rightarrow H^{*-2\mu_d}(S^n M),
\end{equation*}
and $\lambda^*$ is surjective iff $\lambda_r^*$ is surjective. The restriction of the map $\lambda_r^*$ to $H^{*-2\nu_d}(J_d(M))$ is induced by the Abel-Jacobi map $S^n M \rightarrow J_n(M)$.
\end{lem}

\begin{proof}
The same methods as   \cite[Sect.\  7]{AtiyahBott83} show that for the critical set $\eta_d$, the following decomposition of the equivariant cohomology holds
\begin{equation*}
H_\mathcal{G}^* (\eta_d) \cong H_{\mathcal{G}_{diag}}^*(\eta_d) \cong H_{G_{diag}}^*(\tilde{\eta}_d^*)
\end{equation*}
where $\mathcal{G}_{diag}$ is the subgroup of gauge transformations that are diagonal with respect to the Harder-Narasimhan filtration, $\eta_d^*$ refers to the subset of critical points that split with respect to a fixed filtration, $G_{diag}$ is the subgroup of constant gauge transformations that are diagonal with respect to the same fixed filtration, and $\tilde{\eta}_d^*$ is the fiber of $\eta_d^* \cong \mathcal{G}_{diag} \times_{G_{diag}} \tilde{\eta}_d^*$. In the rank $2$ case, the group $G_{diag}$ is simply the torus $T=U(1) \times U(1)$ and we can define (using the local coordinates on $\nu_d^-$ from Section \ref{sec:Kuranishi})
\begin{align}
\tilde{Z}_d^* & = \{ (A, \Phi, a, \varphi) \in (\nu_d^-)_r \, : (A, \Phi) \in \tilde{\eta}_d^*, a = 0 \} \\
Z_d^* & = \{ (A, \Phi, a, \varphi) \in ( \nu_d^-)_r \, : (A, \Phi) \in \tilde{\eta}_d^*, a = 0, \varphi \neq 0 \}
\end{align}
(we henceforth omit the subscript 21 from $(a,\varphi)$; also, $L$ will denote a general line bundle, and not necessarily $L_1\otimes L_2^\ast$).
The map $\lambda^*$ is induced by the inclusion $Z_d^* \hookrightarrow \tilde{Z}_d^*$ and so the map $\lambda^*$ becomes
$
\lambda^* : H_T^*(\tilde{Z}_d^*) \rightarrow H_T^*(Z_d^*)
$.
Let $T'$ be the quotient of $T$ by the subgroup of constant multiples of the identity. Since the constant multiples of the identity fix all points in $\tilde{Z}_d^*$ and $Z_d^*$ then $H_T^*(\tilde{Z}_d^*) \cong H_{T'}^*(\tilde{Z}_d^*) \otimes H^*(BU(1))$ and $H_T^*(Z_d^*) \cong H_{T'}^*(Z_d^*) \otimes H^*(BU(1))$. Therefore the map 
$$\lambda^*: H_{T'}^*(\tilde{Z}_d^*) \otimes H^*(BU(1)) \rightarrow H_{T'}^*(Z_d^*) \otimes H^*(BU(1))$$
 is the identity on the factor $H^*(BU(1))$.

Now consider coordinates on $\tilde{Z}_d^*$ given by $(L_1, L_2, \Phi_1, \Phi_2, \varphi)$ where $L_1\in J_d(M)$,  $L_2 
\in J_{d_E-d}(M)$ are the line bundles of the holomorphic splitting $E = L_1 \oplus L_2$ and $\varphi \in H^0(L_1\otimes L_2^\ast \otimes K)$. For a fixed holomorphic structure, $\Phi_1$ and $\Phi_2$ take values in a vector space, and 
so $\tilde{Z}_d^*$ is homotopy equivalent to
a fibration over 
\begin{equation} \label{eqn:fibration}
\left\{ (L,\varphi) : L\in J_n\ ,\ \varphi\in H^0(L)\right\}
\end{equation}
with fiber $J_d(M)$.  
The fibration is trivialized by  the map
$$
(L_1, L, \varphi) \mapsto (L_1, L_2=L_1\otimes K^\ast\otimes L,\varphi)
$$
 Let $F_n$ be the subspace of \eqref{eqn:fibration} with $\Vert\varphi \Vert=1$.
 Then the cohomology of the fiber bundle splits as
\begin{align}
H_{T'}^*(\tilde{Z}_d^*) & \cong H^*(J_d(M)) \otimes H_{T'}^*(J_n(M)) \label{eqn:etaprimecohomology} \\
H_{T'}^*(Z_d^*) & \cong H^*(J_d(M)) \otimes H_{T'}^*(F_n) \label{eqn:etadoubleprimecohomology}
\end{align}
Note that $F_n$ fibers over the symmetric product $S^n M$ with fiber 
$U(1) \cong T'$, where $T'$ acts trivally on the base, and freely on the fibers. The map $\lambda^*$ restricts 
to the identity on the factor $H^*(J_d(M))$ in \eqref{eqn:etaprimecohomology} and 
\eqref{eqn:etadoubleprimecohomology}, and therefore it restricts to a map $H_{T'}^*(J_n(M)) \rightarrow 
H_{T'}^*(F_n)$. Now the action of $T'$ fixes the holomorphic structures on $L_1$ and $L_2$, and so acts 
trivially on the base of the fiber bundle. $T'$ acts freely on a nonzero section $\varphi \in H^0(L_1^* L_2 
\otimes K)$ and so (after applying the deformation retraction $\left| \varphi \right| \rightarrow 1$), the quotient 
of the space $F_n$ is the space of effective divisors on $M$, since the zeros of each $0 \neq \varphi \in 
H^*(L_1^* L_2 \otimes K)$ correspond to an effective divisor of degree $n = 2g-2 + d_E - 2d$. Therefore the 
map $\lambda^*$ restricts to a map
\begin{equation*}
\lambda_r^* : H^*(J_n(M)) \otimes H^*(BU(1)) \rightarrow H^*(S^n M)
\end{equation*}
which is induced by the $T'$-equivariant map $F_n \rightarrow J_n(M)$, which maps a nonzero section $\varphi \in H^0(L_1^* L_2 \otimes K)$ to the line bundle $L_1^* L_2 \otimes K$. On the quotient $F_n/T' = S^n M$ this restricts to the Abel-Jacobi map $S^n M \rightarrow J_n(M)$.
\end{proof}

Let 
\begin{align*}
\mathcal{M}^{pairs} &= \left\{ (L, \Phi) \, : \, L \in J_n(M),  \Phi \in H^0(L \otimes K) \right\} \\
 \mathcal{M}_0^{pairs} &= \left\{ (L, \Phi) \, : \, L \in J_n(M),  \Phi \in H^0(L \otimes K) \setminus \{0\} \right\}
 \end{align*}
 The group $U(1)$ acts on $\mathcal{M}^{pairs}$ and $\mathcal{M}_0^{pairs}$ by $e^{i\theta} \cdot (L, \Phi) = (L, e^{i\theta} \Phi)$. The inclusion $\mathcal{M}_0^{pairs} \hookrightarrow \mathcal{M}^{pairs}$ is $U(1)$-equivariant with respect to this action, and the proof of Lemma \ref{cla:restrictedsurjective} shows that $\lambda_r^*$ is induced by this inclusion.

\begin{rek}\label{rem:MacDonald}
The paper \cite{MacDonald62} describes the cohomology ring of the symmetric product of a curve  in detail. The result relevant to this paper is that $H^* (S^n M)$ is generated by $2g$ generators in $H^1$, and one generator in $H^2$. Therefore, the proof of Theorem \ref{thm:kirwansurjectivity} reduces to showing that $\lambda_r^*$ maps onto these generators. 
\end{rek}

From the proof of \cite[$\P$ (14.1)]{MacDonald62} we have the following lemma for the Abel-Jacobi map.

\begin{lem}\label{lem:H1surjective}

$\lambda_r^*$ is surjective onto $H^1(S^n M)$.

\end{lem}

Next we need the following technical lemma.

\begin{lem}\label{lem:noirreducible}
For any positive integer $n$, the cohomology group $H^2(F_n)$ consists of products of elements of $H^1(F_n)$.
\end{lem}

\begin{proof}

First consider the case where $n > 2g-2$. By Serre duality $h^1(L) = 0$ for all $L \in J_n(M)$, and so Riemann-Roch shows that $h^0(L) = n+1-g$. Therefore $F_n$ is a sphere bundle over the Jacobian $J_n(M)$ with fiber the sphere $S^{2(n-g+1) - 1}$. By the spectral sequence for this fiber bundle, $H^k(F_n) \cong H^k(J_n(M))$ for all $k \leq 2(n-g+1) - 1$, therefore in low dimensions the ring structure of $H^*(F_n)$ is isomorphic to that of $H^*(J_n(M))$. In particular, since $2(n-g+1)-1 \geq 2g-1 > 2$, we see that $H^2(F_n)$ consists of products of elements of $H^1(F_n)$.

When $n < 2g-2$ we see that $F_n$ is not a fiber bundle over the Jacobian (since the dimension of the fiber may jump). For a fixed basepoint $x_0$ of $M$, consider the inclusion map $M^n \hookrightarrow M^N$ given by 
$$(x_1, \ldots, x_n) \mapsto (x_1, \ldots, x_n, x_0, \ldots, x_0)$$
This induces the inclusion of symmetric products $i: S^n M \hookrightarrow S^N M$, and the description of the generators of $H^*(S^N M)$ in \cite[eq.\ (3.1)]{MacDonald62} shows that the induced map $i^* : H^*(S^N M) \rightarrow H^*(S^n M)$ maps generators to generators and hence is surjective. Therefore the inclusion $i$ induces the following map of fiber bundles
 \begin{equation*}
\left( \begin{matrix} U(1) & \rightarrow & F_n   \\  &  & \downarrow \\  &  & S^n M \end{matrix} \right) \rightarrow \left( \begin{matrix} U(1) & \rightarrow & F_N  \\  &  & \downarrow \\  &  & S^N M \end{matrix} \right)
\end{equation*}
which is the identity map $j: U(1) \rightarrow U(1)$ on the fibers. 

If $N> 2g-2$ then the previous argument implies $H^2(F_N)$ has no irreducible generators, and so in the Serre spectral sequence for $H^*(F_N)$, the irreducible generator $p_N \in H^2 \left(S^N M; H^0(U(1)) \right) \cong H^2 ( S^N M) \otimes H^0(U(1))$ must be killed by a differential (note that $\pi_1(S^N M)$ acts trivially on the space of components of the fiber, and hence on $H^0(U(1))$). For dimensional reasons this must be the differential 
$$
d_2^N : E_2^{0,1}  \cong H^1(U(1)) \otimes H^0(S^N M) \rightarrow E_2^{2,0} \cong H^0(U(1)) \otimes H^2(S^N M)
$$
on the $E_2$ page of the spectral sequence. Since the map $i^*$ is surjective,  $i^* \circ d_2^N$ maps onto $p_n$, the irreducible generator of $H^2(S^n M)$. 

Naturality of the Serre spectral sequence then shows that $d_2^n \circ j^*$ maps onto $p_n$, where $d_2^n : E_2^{0,1} \rightarrow E_2^{2,0}$ is a differential on the $E_2$ page of the Serre spectral sequence for $F_n$. Since $j^*$ is an isomorphism, $d_2^n$ maps onto the irreducible generator $p_n$ of $H^2 \left(S^n M; H^0(U(1)) \right)$.

The following diagram summarizes the argument
\begin{equation*}
\xymatrix{
H^1(U(1)) \otimes H^0(S^N M) \ar[d]^{ j^* \, \mathrm{iso.}} \ar[r]^{d_2^N} & H^0(U(1)) \otimes H^2(S^N M) \ar[d]^{i^* \, \mathrm{surj.}} \\
H^1(U(1)) \otimes H^0(S^n M) \ar[r]^{d_2^n} & H^0(U(1)) \otimes H^2(S^n M)
}
\end{equation*}

Therefore the irreducible generator in $H^2(S^n M)$ is killed by a differential in the spectral sequence for $F_n$, and so there are no irreducible generators of $H^2(F_n)$. 
\end{proof}

\begin{lem}\label{lem:H2surjective}

$\lambda_r^*$ is surjective onto $H^2(S^n M)$.

\end{lem}

\begin{proof}

Using the definition of $F_n$ from above, note that $S^n M \simeq F_n \times_{U(1)} EU(1)$, where $U(1)$ acts by multiplication on the fibers of $U(1) \rightarrow F_n \rightarrow S^n M$. Therefore $S^n M$ is homotopy equivalent to a fiber bundle over $F_n$ with fibers $BU(1)$. From the Serre spectral sequence, we have the map
\begin{multline*}
\left( H^0(F_n) \otimes H^2(BU(1)) \right)  \oplus \left( H^1(F_n) \otimes H^1(BU(1)) \right) \\
 \oplus \left( H^2(F_n) \otimes H^0(BU(1)) \right)  \rightarrow H^2(S^n M)
\end{multline*}
From \cite{MacDonald62}, $H^2(S^n M)$ has an irreducible generator $p_n$. We have that $H^1(BU(1)) = 0$ and by Lemma \ref{lem:noirreducible} there are no irreducible generators of $H^2(F_n) \otimes H^0(BU(1))$. Therefore $p_n$ is in the image of the term $H^0(F_n) \otimes H^2(BU(1)) \cong \C$, and therefore this term is not killed by any differential in the Serre spectral sequence for $S^n M \simeq F_n \times_{U(1)} EU(1)$.

By construction, the map $\lambda_r^*$ is induced by a map of fiber bundles which is an isomorphism on the base $BU(1)$
\begin{equation*}
\left( \begin{matrix}  F_n & \rightarrow & F_n \times_{U(1)} EU(1) \simeq S^n M \\  &  & \downarrow \\  &  & BU(1) \end{matrix} \right) \rightarrow \left( \begin{matrix}  J_n(M) & \rightarrow & J_n(M) \times_{U(1)} EU(1)  \\  &  & \downarrow \\  &  & BU(1) \end{matrix} \right)
\end{equation*}
and therefore the induced map 
$$H^2(BU(1)) \otimes H^0(J_n(M)) \rightarrow H^2(BU(1)) \otimes H^0(F_n)$$
is an isomorphism on the $E_2$ page of the respective Serre spectral sequences. Therefore the map
\begin{equation*}
H^2(BU(1)) \otimes H^0(J_n(M)) \hookrightarrow H^2(J_n(M) \times_{U(1)} EU(1) ) \rightarrow H^2(S^n M)
\end{equation*}
is surjective onto the generator $p_n$ of $H^2(S^nM)$. 
\end{proof}

\begin{proof}[Proof of Theorem \ref{thm:kirwansurjectivity}]
The results of Lemmas \ref{lem:H1surjective} and \ref{lem:H2surjective}, together with MacDonald's results about the cohomology of the symmetric product $S^n M$ (see Remark \ref{rem:MacDonald}), show that the map $\lambda^*$ is surjective. Therefore, Corollary \ref{cor:surjective} implies $\kappa_H$ is surjective.
\end{proof}

\subsection{The action of $\Gamma_2$ on the cohomology}\label{sec:fixedsurjectivity}

First we recall the definition of the action of 
$$\Gamma_2 \cong H^1(M, \mathbb{Z}_2) \cong \Hom(\pi_1(M), \mathbb{Z}_2)$$
 on the space of Higgs bundles (cf.\ \cite{AtiyahBott83, Hitchin87}). $\Gamma_2$ can be identified with the $2$-torsion points of the Jacobian $J_0(M)$ which act on $\mathcal{M}^{Higgs}(2, d_E)$ by tensor product 
 $$L\cdot(E, \Phi) = (E \otimes L, \Phi)$$
  The Jacobian acts also on $\mathcal{M}^{Higgs}(1, k)$ by 
  $$L \cdot (F, \Phi) = (F \otimes L^2, \Phi)$$
   and the determinant map 
$$
\det : \mathcal{M}^{Higgs}(2, d_E)  \rightarrow \mathcal{M}^{Higgs}(1, d_E) :
(E, \Phi)  \mapsto (\det E, \tr \Phi)
$$
becomes $J_0(M)$-equivariant. Since $L \in J_0(M)$ acts on the base by tensoring with $L^2$ we obtain, after lifting $\det$ from $\mathcal{M}^{Higgs}(1, d_E)$ (which is homotopy equivalent to ${J}_0(M)$) to the cover $\widehat {\mathcal{M}}^{Higgs}(1, d_E)$ corresponding to $\Gamma_2$, a product fibration
\begin{equation}\label{eqn:lifted-determinant}
\widehat{\det} : \mathcal{M}_0^{Higgs}(2, d_E) \times \widehat {\mathcal{M}}^{Higgs}(1, d_E) \rightarrow \widehat {\mathcal{M}}^{Higgs}(1, d_E).
\end{equation}
The trivialization
\begin{equation}\label{trivial'}
\hat{\chi} : \mathcal{M}_0^{Higgs}(2, d_E) \times \widehat {\mathcal{M}}^{Higgs}(1, d_E) \rightarrow \widehat {\mathcal{M}}^{Higgs}(2, d_E)
\end{equation}
given by $(E, L)\mapsto E \otimes L$
descends to a homeomorphism
\begin{equation*}
\mathcal{M}_0^{Higgs}(2, d_E) \times_{\Gamma_2} \widehat {\mathcal{M}}^{Higgs}(1, d_E)\cong \mathcal{M}^{Higgs}(2, d_E).
\end{equation*}
(cf.\ \cite[eq.\ (9.5)]{AtiyahBott83} for the case of holomorphic bundles). It is originally one of the main observations of Atiyah and Bott (cf.\ \cite[Sects.\  2 and 9]{AtiyahBott83}) that we can also define the $\Gamma_2$-action via equivariant cohomology. 

Recall from \cite{AtiyahBott83} that the group $\Gamma$ of components of $\mathcal{G}$ is given by $\Gamma \cong H^1(M, \Z)$. Let $\Gamma' = 2\Gamma \subset \Gamma$ be a sublattice of index $2$, and let $\mathcal{G}'$ be the associated subgroup of $\mathcal{G}$, whose components correspond to elements of $\Gamma'$. By  \cite[Prop.\ 2.16]{AtiyahBott83}, $B \mathcal{G}'$ is torsion-free and has the same Poincar\'e polynomial as $B \mathcal{G}$. 

The \emph{degree} of a gauge transformation is the component of $\mathcal{G}$ containing $g$, i.e. $\deg g \in \Gamma$. Dividing by the subgroup of constant central gauge transformations, we obtain $\bar{\mathcal{G}} = \mathcal{G} / U(1)$, and $\bar{\mathcal{G}}_0 = \mathcal{G}_0 / \{ \pm 1\}$, and we define
\begin{equation*}
\bar{\mathcal{G}}' = \{ g \in \bar{\mathcal{G}} \, : \, \deg g \in \Gamma' \} .
\end{equation*}

Let  $\B(1,k)$ denote the space of  Higgs bundles on a line bundle $L\to M$ of degree $k$, $\G(1)$ the corresponding gauge group, and $\G_p(1)$ the subgroup based at $p$.  Fix a basepoint $D_0 \in \mathcal{B}_0(2,d_E)$ and define $T:\B(2,d_E)\to \B(1,d_E)$, the \emph{trace map}, by $T(A,\Phi)=(\tr A, \tr\Phi)$.  Clearly, $T$ is a fibration with fiber $\simeq\B_0(2,d_E)$.  

The fixed determinant gauge group $\G_0$ acts on $\B(2,d_E)$ preserving $\B_0(2,d_E)$ and such that $T$ is invariant.  To see this, note that if
 $g\in \G_0$, then $\tr (D_0 g g^{-1})=0$.  Indeed,
since $\G_0$ is connected it suffices to show that $\tr (D_0 g g^{-1}) = \tr (d g g^{-1})$ is locally constant.  Any $g$ in a  neighborhood of $g_0$ can be expressed $e^u g_0$, where $u\in \Lie(\G_0)$ is a smooth map from $M$ to the vector space of traceless endomorphisms.  In particular, $\tr(du)=d\tr u=0$.  But then 
\begin{align*}
\tr(dg g^{-1})&= \tr(d(e^u)e^{-u})+ \tr(e^u dg_0 g_0^{-1} e^{-u}) \\ 
&=\tr(du)+ \tr( dg_0 g_0^{-1} )= \tr(dg_0 g_0^{-1}) .
\end{align*}
Now for $g\in \G_0$, 
$$
T(g(A), g\Phi g^{-1})=\left( \tr(gAg^{-1} - dg g^{-1}), \tr g\Phi g^{-1}\right) = (\tr A,\tr\Phi) ,
$$
hence there is an induced fibration $T:\B(2,d_E)\times_{\G_0} E\G \to \B(1,d_E)$ with fiber 
$\B_0(2,d_E)\times_{\G_0} E\G$.

The group $\G/\G_0\simeq \G(1)$ induced by the determinant map acts fiberwise on  $T$ with nontrivial stabilizers on $\B(1,d_E)$ given by the constant $\U(1)$ gauge transformations. Therefore, following the approach of  \cite{AtiyahBott83}, we pass to the quotient $\overline\G=\G/\U(1)$, $\overline\G_0=\G_0/\{\pm 1\}$ and consider the induced fibration $\overline T:\B(2,d_E)\times_{\overline\G_0} E\overline\G \to \B(1,d_E)$.  We claim that $\overline T$ is a trivial fibration.  Indeed, with respect to the fixed base point $D_0 \in \B_0(2,d_E)$  define
\begin{align*}
&\chi : \B(2,d_E)\longrightarrow \B_0(2,d_E)\times \B(1,d_E) \\
&\chi(A,\Phi)=\left( (A-(\tfrac{1}{2}\tr A)I, \Phi-(\tfrac{1}{2}\tr\Phi)I), (\tr A, \tr\Phi)\right).
\end{align*}
Then $\chi$ descends to a trivialization
$$
\overline \chi :
\B(2,d_E)\times_{\overline \G_0} E\overline\G\longrightarrow \left(\B_0(2, d_E)\times_{\overline\G_0} E\overline \G\right)\times \B(1,d_E).
$$
Now the group $\overline\G\bigr/\overline\G_0\simeq\overline\G(1)=\G(1)\bigr/\U(1)$ induced by the determinant map acts freely on the total space and the base of $\overline T$, but the induced fibration on the quotient is not trivial.  For this reason we need to pass to a subgroup.

Indeed, given $g\in \G(1)$ let $\deg g\in\Gamma=H^1(M,\mathbb{Z})$ denote the degree of the gauge transformation $g$.  Since constant gauge transformations have degree $0$, it induces a map $\deg: \overline \G(1)\to \Gamma$.  Let
$$
\overline \G'(1)=\left\{ g\in\overline \G(1) : \deg g\in 2\Gamma\right\} .
$$
We  define  $\overline \G'= \left\{g\in \overline \G : \det(g) \in \overline \G'(1)\right\}$.  
Given $g\in \overline\G'(1)$, set $g=s^2$, $s\in \overline \G(1)$, and let $\hat g=\left(\begin{matrix} s&0\\ 0&s\end{matrix}\right)\in \overline \G$.  Define $g[A,\Phi, e]=[\hat g(A,\Phi, e)]$ for  $[A,\Phi, e]\in\B(2,d_E)\times_{\overline\G_0} E\overline\G$.  Notice that the action is well-defined independent of the choice of square root.  Furthermore, $\overline \chi$ is equivariant, where the action of $\overline \G'(1)$ is trivial on $\B_0(2,d_E)\times_{\overline\G_0} E\overline\G$ and has the usual action on $\B(1,d_E)$.  Hence the induced fibration
\begin{equation}\label{trivialfibr}
\widehat T : \B(2,d_E)\times_{\overline\G'} E\overline \G \longrightarrow \B(1,d_E)\bigr/\overline\G'(1)
\end{equation}
can be trivialized by the homeomorphism
\begin{equation}\label{trivializ}
\hat \chi :\B(2,d_E)\times_{\overline \G'} E\overline\G\longrightarrow \left(\B_0(2, d_E)\times_{\overline\G_0} E\overline \G\right)\times \B(1,d_E)\bigr/\overline\G'(1)
\end{equation}
induced from $\overline \chi$. 
\begin{rek}
Formulas \eqref{trivialfibr} and \eqref{trivializ} should be considered as the equivariant analogues of \eqref{eqn:lifted-determinant} and \eqref{trivial'}.
\end{rek}

Now $\Gamma_2$ acts on the left hand side of (\ref{trivializ}). It is also clear that the action of $\Gamma_2$ on  $B(1,d_E)\bigr/\overline\G'(1)\cong\widehat {\mathcal{M}}^{Higgs}(1, d_E)$
is just by tensoring with a torsion point in the Jacobian. 
\begin{defi}\label{defact} The action of $\Gamma_2$ on  $\B_0(2,d_E)\times_{\overline\G_0} E\overline \G$ is defined 
so that the map $\hat \chi$ becomes $\Gamma_2$-equivariant.
\end{defi}
The following simple lemma identifies also the two actions on the fibers of
\eqref{eqn:lifted-determinant} and \eqref{trivialfibr}.
\begin{lem}\label{act} On any subspace $Y$ of $\B_0(2,d_E)$ invariant under $\G_0$ on which $\G_0$ acts with constant stabilizer, the action of $\Gamma_2$ 
on $Y/\G_0$ is given by
tensoring with a 2-torsion point of $J_0(M)$.
\end{lem}
\begin{proof}
Given $\gamma \in \Gamma_2$, let $g_{\gamma}$ be a gauge transformation in $\overline \G(1)$ such that $\deg(g_{\gamma})=\gamma \mod  H^1(M, 2\mathbb Z) $ and 
$h_{\gamma}\in\overline \G$ with $\det (h_{\gamma})= g_{\gamma}$. Note that $\tr (h_\gamma^{-1} d h_\gamma) = g_\gamma^{-1} d g_\gamma$. Then by Definition \ref{defact}, the action of $h_\gamma$ on $\mathcal{B}_0(2, d_E)$ (modulo gauge transformations in $\mathcal{G}_0$) is given by
\begin{align*}
h_\gamma \left[ (A, \Phi ) \right] & = [  ( h_\gamma^{-1} A h_\gamma + h_\gamma^{-1} D_0 h_\gamma - \frac{1}{2} \tr \left( h_\gamma^{-1} D_0 h_\gamma + h_\gamma^{-1} A h_\gamma \right) I, \\
 & \quad \quad h_\gamma^{-1} \Phi h_\gamma - \frac{1}{2} \tr \left( h_\gamma^{-1} \Phi h_\gamma \right) I)] \\
 & = \left[ \left( h_\gamma^{-1} D_0 h_\gamma + h_\gamma^{-1} A h_\gamma - \frac{1}{2} (g_\gamma^{-1} d g_\gamma) I , h_\gamma^{-1} \Phi h_\gamma \right) \right],
\end{align*}
since $\tr A = 0$ and $\tr \Phi = 0$. We claim that this equivalent to tensoring with the line bundle $L_{\gamma}$ corresponding to $\gamma$. To see this last statement, chose a simple loop $\sigma$ on $M$ and note that if $\gamma[\sigma]=+1$, then $g_\gamma$ has even degree around the loop $\sigma$ and so in an annulus around $\sigma$ the gauge transformation $g_{\gamma}=s^2$ is a square, hence the previous formula becomes
$$
h_\gamma [ (A, \Phi) ] = \left[ (\hat{g}^{-1}  h_\gamma) \cdot (A, \phi) \right] ,
$$ 
where $\hat{g} = s I$ as before (note that since $g_\gamma \in \bar{\mathcal{G}}(1)$ then $g_\gamma^{-1} d g_\gamma = d g_\gamma g^{-1}$).  Since $\hat{g}^{-1} h_\gamma \in \mathcal{G}_0$ then this shows that $h_\gamma [(A, \Phi)] = [(A, \Phi)]$ in an annulus around $\sigma$.
 
If $\gamma[\sigma]=-1$ then parametrise the loop $\sigma$ by $\theta \, : 0 \leq \theta \leq 2 \pi$ and note that since $g_\gamma$ has odd degree, then $g_\gamma = e^{i \theta} s^2$ in an annulus around $\sigma$. Therefore the effect of the gauge term  $(\frac{1}{2} {g_{\gamma}}^{-1}d{g_{\gamma}})I$ is that it changes the argument of the holonomy around $\sigma$ by $\pi$, as desired.
\end{proof}

In the above we can restrict to the $\mathcal{G}_0$-invariant subspaces $X_d$ of $\mathcal{B}_0(2, d_E)$, and the action commutes with inclusions and connecting homomorphisms from the LES in cohomology. Therefore, we have a LES of $\Gamma_2$ spaces and $\Gamma_2$-equivariant maps
\begin{equation*}
\xymatrix{
 H_{\mathcal{G}_0}^k(X_d, X_{d-1}) \ar[r]^(0.6){\alpha^k} & H_{\mathcal{G}_0}^k(X_d) \ar[r]^{\beta^k} & H_{\mathcal{G}_0}^k(X_{d-1}) \ar[r]^(0.4){\gamma^k}  & H_{\mathcal{G}_0}^{k+1} (X_d, X_{d-1}) 
}
\end{equation*}

\begin{lem}\label{lem:gamma2-action-cohomology}
The $\Gamma_2$-action commutes with the isomorphism in \eqref{eqn:bott}
\begin{equation}\label{eqn:excision}
H_{\mathcal{G}_0}^* (X_d, X_{d-1}) \cong H_{\mathcal{G}_0}^*(\nu_d^-, \nu_d'),
\end{equation}
and with the isomorphisms  \eqref{eqn:nucohomology-doubleprime} and \eqref{eqn:nucohomology-prime-doubleprime}, \eqref{eqn:fixedcohomologycritical} and \eqref{eqn:fixednucohomology-prime-doubleprime}.
\end{lem} 
\begin{proof}
First, note that the $\Gamma_2$ action on $\mathcal{B}_0(2, d_E) \times_{\bar{\mathcal{G}_0}} E \bar{\mathcal{G}_0}$ preserves the subspaces $\mathcal{B}_d \times_{\bar{\mathcal{G}_0}} E\bar{\mathcal{G}_0}$ and $\nu_d^- \times_{\bar{\mathcal{G}_0}} E \bar{\mathcal{G}_0}$, $X_d \times_{\bar{\mathcal{G}_0}} E \bar{\mathcal{G}_0}$ and $\nu_d' \times_{\bar{\mathcal{G}_0}} E \bar{\mathcal{G}_0}$ for all values of $d$, and so the inclusion of pairs
\begin{equation*}
\left(\nu_d^- \times_{\bar{\mathcal{G}_0}} E \bar{\mathcal{G}_0}, \nu_d' \times_{\bar{\mathcal{G}_0}} E \bar{\mathcal{G}_0} \right) \hookrightarrow \left( X_d \times_{\bar{\mathcal{G}_0}} E \bar{\mathcal{G}_0}, X_{d-1} \times_{\bar{\mathcal{G}_0}} E \bar{\mathcal{G}_0} \right)
\end{equation*}
is $\Gamma_2$-equivariant. Therefore the action of $\Gamma_2$ commutes with the excision isomorphism
\begin{equation*}
H^*(X_d \times_{\bar{\mathcal{G}_0}} E \bar{\mathcal{G}_0}, X_{d-1} \times_{\bar{\mathcal{G}_0}} E \bar{\mathcal{G}_0}) \cong H^*(\nu_d^- \times_{\bar{\mathcal{G}_0}} E \bar{\mathcal{G}_0}, \nu_d' \times_{\bar{\mathcal{G}_0}} E \bar{\mathcal{G}_0})
\end{equation*}
which descends to the isomorphism \eqref{eqn:excision} in equivariant cohomology.

The isomorphisms \eqref{eqn:fixedcohomologycritical} and \eqref{eqn:fixednucohomology-prime-doubleprime} arise from taking quotients 
\begin{equation*}\label{eqn:fixed-critical-gamma}
H_{\bar{\mathcal{G}_0}}^*(\eta_d) \cong H^*(\eta_d \times_{\bar{\mathcal{G}_0}} E \bar{\mathcal{G}_0}) \cong H_{U(1)}^* ( \eta_d / \bar{\mathcal{G}_0}) \cong H^* \left(J_d(M)\right) \otimes H^* \left(BU(1) \right)
\end{equation*}
(where $\bar{\mathcal{G}_0}$ acts on $\eta_d$ with isotropy group $\U(1)$), and 
\begin{equation}\label{eqn:fixed-nonzero-bad-gamma}
H_{\bar{\mathcal{G}_0}}^*(T_d) \cong H^*( T_d \times_{\bar{\mathcal{G}_0}} E \bar{\mathcal{G}_0} )\cong H^*(T_d / \bar{\mathcal{G}_0}) \cong H^* ( \widetilde{S}^n M)
\end{equation}
(since $\bar{\mathcal{G}_0}$ acts freely on $T_d$). The action of $\Gamma_2$ on the space $\mathcal{B}_0 \times_{\bar{\mathcal{G}_0}} E \bar{\mathcal{G}_0}$ induces actions on $\eta_d \times_{\bar{\mathcal{G}_0}} E \bar{\mathcal{G}_0}$ and $T_d \times_{\bar{\mathcal{G}_0}} E \bar{\mathcal{G}_0}$ which in turn induces an action on the spaces $\eta_d / \bar{\mathcal{G}_0}$ and $T_d / \bar{\mathcal{G}_0}$. 
By Lemma \ref{act} the action of $\gamma \in \Gamma_2$ on the quotient $\eta_d / \bar{\mathcal{G}_0} \simeq \left\{ (L_1, L_2) \in J_d(M) \times J_{d_E-d}(M) \, : \, L_1 L_2 = F \right\}$,  is given by tensor product $(L_1, L_2) \mapsto (L_1 \otimes L_{\gamma}, L_2 \otimes L_{\gamma})$, where $L_{\gamma} \in J_0(M)$ is the line bundle corresponding to $\gamma$. The induced action on the cohomology is trivial by  \cite[Prop.\ 9.7]{AtiyahBott83}. The action of $\Gamma_2$ on the quotient $T_d / \bar{\mathcal{G}_0}$ is also by tensor product, $(L_1, L_2, \Phi) \mapsto (L_1 \otimes L_{\gamma}, L_2 \otimes L_{\gamma}, \Phi)$, therefore the action on the right-hand side of \eqref{eqn:fixed-nonzero-bad-gamma} is via deck transformations of the $2^{2g}$-fold cover $\widetilde{S}^n M \rightarrow S^n M$ (see also \cite[Sect.\  7]{Hitchin87}).
\end{proof}

Let $N$ be a space with a $\Gamma_2$-action.  Then we have a  splitting $$H^*(N) \cong H^*(N)^{\Gamma_2} \oplus H^*(N)^a$$
where $H^*(N)^{\Gamma_2}$ is the $\Gamma_2$-invariant part of the cohomology and 
$$H^*(N)^a \cong \oplus_{\varphi \neq 1} H^*(N)_{\varphi}$$ where $\varphi$  varies over all homomorphisms $\Gamma_2\to \{\pm 1\}$.   If $N_1$, $N_2$ are two such spaces and 
$f: H^\ast(N_1)\to H^\ast(N_2)$ is a
$\Gamma_2$-equivariant homomorphism, we denote by $f_{\Gamma_2}$ (resp.\ $f_a$) the restriction of $f$ to $H^*(N_1)^{\Gamma_2}$ (resp.\ $H^*(N_1)^a$).

Applying this notation to $\lambda^*$ we have
\begin{equation*}
\lambda_{\Gamma_2}^* : H_\mathcal{G}^*(\nu_d^-, \nu_d'')^{\Gamma_2} \rightarrow H_\mathcal{G}^*(\nu_d', \nu_d'')^{\Gamma_2} .
\end{equation*}
The main result of this section is Lemma \ref{lem:lambda-invariant} which shows that $\lambda_{\Gamma_2}^*$ is surjective, a key step towards proving Theorem \ref{thm:invariantkirwansurjectivity}. The earlier results \eqref{eqn:nucohomology-prime-doubleprime} and Lemma \ref{lem:critsetcohomology} show that $H_\mathcal{G}^*(\nu_d', \nu_d'') \cong H^{*-2\mu_d}(\widetilde{S}^n M)$, where $n=2g-2+d_E-2d$. Points in $\widetilde{S}^n M$ correspond to triples $(L_1, L_2, \Phi) \in J_d(M) \times J_{d_E-d}(M) \times \Omega^0(L_1^* L_2 \otimes K)$ where $L_1 L_2 = \det E$ is a fixed line bundle. Similarly, there is a corresponding $2^{2g}$ cover of the Jacobian $\tilde{J}_n(M) = J_d(M) \times J_{d_E-d}(M)/\negthickspace\sim$, where the equivalence is given by $(L_1, L_2) \sim (\tilde{L}_1, \tilde{L}_2)$ if $L_1 L_2 \cong \tilde{L}_1 \tilde{L}_2$.  

The isomorphisms 
\begin{align*}
H_\mathcal{G}^*(\nu_d^-, \nu_d'') & \cong H^{* - 2\mu_d}(\eta_d) \cong H^{*-2\mu_d}( \tilde{J}(M) \times BU(1) ) \\
H_\mathcal{G}^*(\nu_d', \nu_d'') & \cong H^{* - 2\mu_d}(T_d) \cong H^{*-2\mu_d}(\widetilde{S}^n M)
\end{align*}
from Theorem \ref{thm:nucohomology} and Lemma \ref{lem:critsetcohomology} show that the map $\lambda_{\Gamma_2}$  is given by
\begin{equation*}
\xymatrix{
\lambda_{\Gamma_2}^* : H^{*-2\mu_d}( \tilde{J}(M) \times BU(1) )^{\Gamma_2} \ar[r] \ar[d]^{\cong} & H^{*-2\mu_d}(\widetilde{S}^n M)^{\Gamma_2} \cong H^{*-2\mu_d}(S^n M) \ar[d]^{\cong} \\
H^{*-2 \mu_d}(\eta_d)^{\Gamma_2} \ar[r] & H^{*-2\mu_d}(T_d)^{\Gamma_2}   }
\end{equation*}
where $n = 2g-2+d_E-2d$, and $\mu_d = 2d-d_E+g-1$. This map is induced by the inclusion $T_d \hookrightarrow (\nu_d^-)_r$ (where the spaces are now subsets of the space of fixed determinant Higgs bundles). We define the \emph{lifted Abel-Jacobi map} to be the map $\widetilde{S}^n M \rightarrow \tilde{J}(M)$, which takes a triple $(L_1, L_2, \Phi)$ to the pair $(L_1, L_2) \in \tilde{J}(M)$. The same proof as Lemma \ref{cla:restrictedsurjective} in the previous section gives us the following 
\begin{lem}\label{lem:lifted-abel-jacobi}
The restriction of  $\lambda_{\Gamma_2}^*$ to $H^{*-2\mu_d}(\tilde{J}_n(M))$ given by 
$$
(\lambda_{\Gamma_2}^*)_r : H^{*-2\mu_d}( \tilde{J}_n(M) )^{\Gamma_2} \rightarrow H^{*-2\mu_d}(\widetilde{S}^n M)^{\Gamma_2}
$$
 is induced by the lifted Abel-Jacobi map.
\end{lem}

\begin{lem}\label{lem:lambda-invariant}
The map $\lambda_{\Gamma_2}^*$ is surjective. 
\end{lem}

\begin{proof}
By \cite[eqs.\ (7.12) and (7.13)]{Hitchin87}, $H^*(\widetilde{S}^n M)^{\Gamma_2} \cong H^*(S^n M)$, and we also have $H^*(\tilde{J}_n(M))^{\Gamma_2} \cong H^*(J_n(M))$. Therefore Lemma \ref{lem:H1surjective} implies $\lambda_{\Gamma_2}^*$ is surjective onto $H^1(\widetilde{S}^n M)^{\Gamma_2}$. By the same argument as in Lemma \ref{lem:H2surjective} (with the $\Gamma_2$-invariant part of the cohomology), $\lambda_{\Gamma_2}^*$ is surjective onto $H^2(\widetilde{S}^n M)^{\Gamma_2}$. By \cite{MacDonald62}, $H^*(\widetilde{S}^n M)^{\Gamma_2} \cong H^*(S^n M)$ is generated in dimensions $1$ and $2$; hence, $\lambda_{\Gamma_2}^*$ is surjective.
\end{proof}

\subsection{$\Gamma_2$-invariant hyperk\"ahler Kirwan surjectivity}
For fixed determinant  the inclusion $\mathcal{B}_0^{ss} \hookrightarrow T^* \mathcal{A}_0$ induces a map on the $\Gamma_2$-invariant part of the $\mathcal{G}$-equivariant cohomology which we call the \emph{$\Gamma_2$-invariant hyperk\"ahler Kirwan map}
\begin{equation*}
\kappa_{HK}^{\Gamma_2} : H_\mathcal{G}^* (T^* \mathcal{A}_0) \cong H_\mathcal{G}^*(T^* \mathcal{A}_0)^{\Gamma_2} \rightarrow H_\mathcal{G}^* (\mathcal{B}_0^{ss} )^{\Gamma_2}.
\end{equation*}
In this section we prove
\begin{thm}\label{thm:invariantkirwansurjectivity}
 $\kappa_{HK}^{\Gamma_2}$ is surjective.
\end{thm}
\noindent As mentioned in the Introduction, it turns out that the full Kirwan map is \emph{not} surjective.

The second goal of this section is the following.  The
 results of Section \ref{sec:surjectivity} show that the map $\zeta^k$ in Diagram \eqref{eqn:diagram} is always injective for non-fixed determinant Higgs bundles, and so Lemma \ref{lem:subset} implies that in this case $\ker \alpha^k \cong \ker \zeta^k = \{ 0 \}$. In this section we will show that $\ker \alpha^k \cong \ker \zeta^k$ holds  for fixed determinant as well, which is important for the calculations in Section \ref{sec:computations}.
\begin{prop}\label{prop:alpha}
For rank $2$ Higgs bundles, $\ker \alpha^k \cong \ker \zeta^k$ for all $k$, and therefore $\dim \, \im \alpha^k = \dim \, \im \zeta^k$ also. In the non-fixed determinant case $\ker \alpha^k = 0$ for all $k$, and in the fixed determinant case 
\begin{align*}
\ker \alpha^k & = H_\mathcal{G}^k(X_d, X_{d-1})^a \\
 & \cong \begin{cases} H^{k-2\mu_d}(\widetilde{S}^{2g-2d-2+d_E} M)^a & k = 4g-4-d_E+2d+1 \\ 0 & \mathrm{otherwise} \end{cases} 
\end{align*}
\end{prop}
\noindent  Note that we have already proven $\ker \alpha^k \cong \ker \zeta^k$ in the non-fixed determinant case (Sect.\ \ref{sec:surjectivity}).  Hence, for the rest of this section  we restrict  to the fixed determinant case.

In order to separate out the $\Gamma_2$-invariant part of the equivariant cohomology, we require the following simple  

\begin{lem} \label{lem:isotypical-exact-sequence}
Let
$$
\xymatrix{\cdots\ar[r] & A_n \ar[r]^{f_n} & B_n \ar[r]^{g_n} & C_n \ar[r]^(0.4){h_n} &  A_{n+1} \ar[r] & \cdots}
$$
be a LES of $\mathbb{C}$-vector spaces.  Suppose that $\Gamma$ is a finite abelian group acting linearly on $A_n$, $B_n$ and $C_n$ such that  $f_n$, $g_n$, and $h_n$ are equivariant. Then for each homomorphism ${\varphi}: \Gamma\to \mathbb{C}^\ast$ the restriction
$$
\xymatrix{\cdots\ar[r] & (A_n)_\varphi \ar[r]^{f_{n,\varphi}} & (B_n)_\varphi \ar[r]^{g_{n,\varphi}} & (C_n)_\varphi \ar[r]^(0.4){h_{n,\varphi}} &  (A_{n+1})_\varphi\ar[r]& \cdots}
$$
to the ${\varphi}$-isotypical subspaces is exact.
\end{lem}
\begin{proof}
By the equivariance of the maps the restrictions are well-defined.  We prove exactness at $(B_n)_\varphi$. By equivariance and exactness of the original sequence,   
$$f_n((A_n)_\varphi)\subset\ker g_n\cap (B_n)_{\varphi}$$
Suppose $ b\in\ker g_n\cap (B_n)_{\varphi}$.  Again by exactness of the original sequence, $b=f_n(\tilde a)$ for some $\tilde a\in A_n$.  Set
$$
a=\frac{1}{\#\Gamma}\sum_{\sigma\in \Gamma} {\varphi}(\sigma^{-1})\sigma\tilde a
$$
Then
$$
f_n(a)=\frac{1}{\#\Gamma}\sum_{\sigma\in \Gamma} {\varphi}(\sigma^{-1})\sigma b
=\frac{1}{\#\Gamma}\sum_{\sigma\in \Gamma}{\varphi}(\sigma^{-1}) {\varphi}(\sigma) b=\frac{1}{\#\Gamma}\sum_{\sigma\in \Gamma} b= b
$$
and  since $b\in (B_n)_{\varphi}$,
$$
\gamma a=\frac{1}{\#\Gamma}\sum_{\sigma\in \Gamma} {\varphi}(\sigma^{-1})\gamma\sigma\tilde a
=
\frac{1}{\#\Gamma}\sum_{\gamma\sigma\in \Gamma} {\varphi}((\gamma\sigma)^{-1}){\varphi}(\gamma)\gamma\sigma\tilde a
={\varphi}(\gamma) a
$$
Hence, $a\in (A_n)_{\varphi}$ and $f_n(a)=b$. This completes the proof. 
\end{proof}

We  apply this result to the vertical and horizontal long exact sequences in \eqref{eqn:diagram}.  

\begin{prop}\label{prop:vertical-LES}
The decomposition of the {\bf vertical} LES of Diagram \eqref{eqn:diagram}  into $\Gamma_2$-invariant and noninvariant parts gives the following for all $k$:
\begin{enumerate}
\item 
 $\delta^k_a: H_\mathcal{G}^{k-1}(\nu_d', \nu_d'')^a\to H_\mathcal{G}^k(\nu_d^-, \nu_d')^a$ is an isomorphism; in particular, $H_{\mathcal{G}}^*(\nu_d^-, \nu_d'')^a = 0$.
\item  the sequence
\begin{equation*}\label{eqn:vertical-invariant-exact}
\xymatrix{
0 \ar[r] & H_\mathcal{G}^k(\nu_d^-, \nu_d')^{\Gamma_2} \ar[r]^{\zeta^k_{\Gamma_2}} & H_\mathcal{G}^k(\nu_d^-, \nu_d'')^{\Gamma_2}  \ar[r]^{\lambda^k_{\Gamma_2}} & H_\mathcal{G}^k(\nu_d', \nu_d'')^{\Gamma_2} \ar[r]^(0.7){\delta^k_{\Gamma_2}} & 0
}
\end{equation*}
is exact.
\end{enumerate}
\end{prop}

\begin{proof}
Since the $\Gamma_2$ action is trivial on the cohomology of the Jacobian and on the cohomology of $BU(1)$, it follows from  \eqref{eqn:nucohomology-doubleprime} and \eqref{eqn:fixedcohomologycritical} that $H_{\mathcal{G}}^*(\nu_d^-, \nu_d'')^a = 0$.
Lemma \ref{lem:lambda-invariant} implies $H_\mathcal{G}^k(\nu_d', \nu_d'')^{\Gamma_2} \subseteq \im \lambda^k = \ker \delta^k$, so $\delta^k_{\Gamma_2}=0$ for all $k$, which proves the second part of the Proposition.
The first part then follows from Lemma \ref{lem:isotypical-exact-sequence}.
\end{proof}

\begin{cor}\label{cor:alpha-invariant-injective}
$\alpha^k_{\Gamma_2}$ is injective. 
\end{cor}

\begin{proof}
Let $x \in H_\mathcal{G}^k(X_d, X_{d-1})^{\Gamma_2} \cong H_\mathcal{G}^k(\nu_d^-, \nu_d')^{\Gamma_2}$, and suppose that $\alpha^k(x) = 0$. 
In the following,  use $x$ to also denote the corresponding element in $H_\mathcal{G}^k(\mathcal{B}_{d, \varepsilon}, \mathcal{B}_{d, \varepsilon}')$ via the excision isomorphism. Then from the commutativity of Diagram \eqref{eqn:diagram},  $\alpha^k(x) = 0$ implies that $\alpha_\varepsilon^k(x) = 0$, and so $\xi^k \circ \zeta^k(x) = 0$. By Lemma \ref{lem:Eulerclass} and Proposition \ref{prop:vertical-LES}, $\xi^k$ is injective and  $\zeta^k$ is injective on $H_\mathcal{G}^k(\nu_d^-, \nu_d')^{\Gamma_2}$. Therefore $x = 0$, which completes the proof.
\end{proof}

Lemma \ref{lem:critsetcohomology}, Theorem \ref{thm:nucohomology}, and Lemma \ref{lem:gamma2-action-cohomology}, together with Hitchin's formulas  \cite[eqs.\ (7.12) and (7.13)] {Hitchin87}, give us the following result. 

\begin{lem}\label{lem:Symmetric-product-invariant-part}
\begin{equation*}
H_{\mathcal{G}}^k (\nu_d', \nu_d'')^a =  \begin{cases} V & k = 4g-4 - d_E + 2d \\ 0 & \mathrm{otherwise} \end{cases} 
\end{equation*}
where $V \cong H^{k-2\mu_d}(\widetilde{S}^{2g-2d-2+d_E} M)^a$ is a complex vector space of dimension $$\dim_\C V = (2^{2g}-1)\left( \begin{matrix} 2g-1 \\ 2g-2d-2+d_E \end{matrix} \right).$$
\end{lem}

\begin{lem}\label{lem:low-degree-vanishing}
$
H_\mathcal{G}^k(X_d)^a = 0
$,
for all $k \leq 4g-4-d_E+2d+1$.
\end{lem}

\begin{proof}
The proof is by induction on the index $d$. For $d > g-1$ the induced map $\kappa_H : H_\mathcal{G}^* (\mathcal{B}) \rightarrow H_\mathcal{G}^*(X_d)$ is surjective, since each stratum has a well-defined normal bundle, and so the methods of \cite{AtiyahBott83} work in this case. Therefore, when $d > g-1$ we have that $H_\mathcal{G}^*(X_d)$ is $\Gamma_2$-invariant for all $k$.
Suppose the result is true for $X_d$. To complete the induction we show that it is true for $X_{d-1}$, i.e. $H_\mathcal{G}^k(X_{d-1})$ is $\Gamma_2$-invariant for all $k \leq 4g-4-d_E+2d-1$.

Consider the following LES for $k \leq 4g-4-d_E+2d-1$.
\begin{equation}\label{eqn:low-degree-sequence}
\xymatrix{
\cdots \ar[r]^{\alpha^k} & H_\mathcal{G}^k (X_d) \ar[r]^{\beta^k} & H_\mathcal{G}^k(X_{d-1}) \ar[r]^(0.4){\gamma^k} & H_\mathcal{G}^{k+1}(X_d, X_{d-1}) \ar[r]^(0.7){\alpha^{k+1}} \ar[r] & \cdots
}
\end{equation}
 From Lemma \ref{lem:Symmetric-product-invariant-part} and Proposition \ref{prop:vertical-LES} we see that
\begin{equation*}
H_\mathcal{G}^{k+1}(X_d, X_{d-1})^a \cong H_\mathcal{G}^{k+1}(\nu_d^-, \nu_d')^a \cong H_\mathcal{G}^{k}(\nu_d', \nu_d'')^a = 0
\end{equation*}
for all $k \leq 4g-4-d_E+2d-1$. Therefore $H_\mathcal{G}^{k+1}(X_d, X_{d-1})$ is $\Gamma_2$-invariant. The exact sequence \eqref{eqn:low-degree-sequence} decomposes to become
\begin{equation*}
\xymatrix{
0 \ar[r] & \im \beta^k \ar[r] & H_\mathcal{G}^k (X_{d-1}) \ar[r]^{\gamma^k} & \im \gamma^k \ar[r] & 0
}
\end{equation*}
Since $\im \gamma^k \subseteq H_\mathcal{G}^{k+1}(X_d, X_{d-1})$, and the latter is $\Gamma_2$-invariant, an application of Lemma \ref{lem:isotypical-exact-sequence} implies
\begin{equation*}
0 \longrightarrow (\im \beta^k)^a \longrightarrow H_\mathcal{G}^k (X_{d-1})^a \longrightarrow 0
\end{equation*}
is exact. By the inductive hypothesis, $H_\mathcal{G}^k(X_d)$ is $\Gamma_2$-invariant; hence, $(\im \beta^k)^a = 0$, and so $H_\mathcal{G}^k (X_{d-1})^a = 0$  also.
\end{proof}

\begin{prop}\label{prop:horizontal-LES}
The decomposition of the {\bf horizontal} LES of Diagram \eqref{eqn:diagram}  into $\Gamma_2$-invariant and noninvariant parts gives the following for all $k \leq 4g-4-d_E+2d+1$:
\begin{enumerate}
\item 
 $\gamma^{k-1}_a: H_\mathcal{G}^{k-1}(X_{d-1})^a\to  H_\mathcal{G}^k(X_d, X_{d-1})^a $ is an isomorphism; in particular, $H_\mathcal{G}^{k-1}(X_{d-1})^a \cong H_\mathcal{G}^{k}(\nu_d^-, \nu_d')^a$.
 \item  the sequence
\begin{equation*}\label{eqn:horizontal-invariant-exact}
\xymatrix{
0 \ar[r] & H_\mathcal{G}^{k-1}(X_d, X_{d-1})^{\Gamma_2} \ar[r]^(0.6){\alpha^{k-1}_{\Gamma_2}} & H_\mathcal{G}^{k-1}(X_d)^{\Gamma_2}  \ar[r]^{\beta^{k-1}_{\Gamma_2}} & H_\mathcal{G}^{k-1}(X_{d-1})^{\Gamma_2} \ar[r]^(0.7){\gamma^{k-1}_{\Gamma_2}} & 0
}
\end{equation*}
is exact.
\end{enumerate}
\end{prop}

\begin{proof}
First, by Lemma \ref{lem:low-degree-vanishing}, $H_\mathcal{G}^{k-1}(X_d)^a = 0 = H_\mathcal{G}^k (X_d)^a$ for $k \leq 4g-4-d_E+2d+1$. 
Next we claim that $\gamma^{k-1}$ maps $H_\mathcal{G}^{k-1}(X_{d-1})^{\Gamma_2}$ to zero for all values of $k$ (not just for $k \leq 4g-4-d_E+2d+1$). To see this, let $x \in H_\mathcal{G}^{k-1}(X_{d-1})^{\Gamma_2}$, and let $y = \gamma^{k-1}(x) \in H_\mathcal{G}^k(X_d, X_{d-1})^{\Gamma_2}$. Exactness of the horizontal LES in Diagram \eqref{eqn:diagram} implies $\alpha^k (y) = \alpha^k \circ \gamma^{k-1}(x) = 0$. By Corollary \ref{cor:alpha-invariant-injective},  $\alpha^k$ is injective on $H_\mathcal{G}^k(X_d, X_{d-1})^{\Gamma_2}$; hence, $y = \gamma^{k-1}(x) = 0$. Therefore, $\gamma^{k-1}(x) = 0$, and so $\gamma^{k-1}$ is the zero map on $H_\mathcal{G}^{k-1}(X_{d-1})^{\Gamma_2}$.
The result then follows from Lemma \ref{lem:isotypical-exact-sequence}.
\end{proof}

\begin{proof}[Proof of Theorem \ref{thm:invariantkirwansurjectivity}]
By the proof of Proposition \ref{prop:horizontal-LES}, $\gamma^k_{\Gamma_2}=0$  for all $k$. By Lemma \ref{lem:Symmetric-product-invariant-part}, $H_\mathcal{G}^k(X_{d-1})^a$ is only nontrivial for $k = 4g-4-d_E+2d$, and so 
Proposition \ref{prop:horizontal-LES} (i) implies $\gamma^k$ is injective on $H_\mathcal{G}^k(X_{d-1})^a$ for all $k$. Therefore, $\beta^k$ maps $H_\mathcal{G}^k(X_d)^{\Gamma_2}$ surjectively onto $H_\mathcal{G}^k(X_{d-1})^{\Gamma_2}$ for all $k$.
Applying this result to every stratum $X_d$ completes the proof of the theorem.
\end{proof}

\begin{proof}[Proof of Proposition \ref{prop:alpha}]
For $k \leq 4g-4-d_E+2d+1$,  Proposition \ref{prop:horizontal-LES} (i) implies $\ker \alpha^k \supseteq H_\mathcal{G}^k(X_d, X_{d-1})^a$, which together with Corollary \ref{cor:alpha-invariant-injective} implies $\ker \alpha^k = H_\mathcal{G}^k(X_d, X_{d-1})^a \cong H_\mathcal{G}^k(\nu_d^-, \nu_d')^a$. The two exact sequences in Proposition \ref{prop:vertical-LES} show that $\ker \zeta^k = H_\mathcal{G}^k(\nu_d^-, \nu_d')^a \cong H_\mathcal{G}^k(\nu_d', \nu_d'')^a$. Therefore Lemma \ref{lem:Symmetric-product-invariant-part} implies
\begin{align*}
\ker \alpha^k \cong \ker \zeta^k & \cong H_\mathcal{G}^{k-1}(\nu_d', \nu_d'')^a \\
 & \cong \left\{ \begin{matrix} H^{k-1-2\mu_d}(\widetilde{S}^{2g-2d-2+d_E} M)^a & k = 4g-4-d_E+2d+1 \\ 0 & k < 4g-4-d_E+2d+1 \end{matrix} \right.
\end{align*}

For $k > 4g-4-d_E+2d+1$, Lemma \ref{lem:Symmetric-product-invariant-part} and Proposition  \ref{prop:vertical-LES} show that 
$$H_\mathcal{G}^k (X_d, X_{d-1})^a \cong H_\mathcal{G}^k(\nu_d^-, \nu_d')^a \cong H_\mathcal{G}^{k-1}(\nu_d', \nu_d'')^a = 0$$ 
Hence,
$H_\mathcal{G}^k(X_d, X_{d-1}) = H_\mathcal{G}^k(X_d, X_{d-1})^{\Gamma_2}$,
and so   $\ker \alpha^k = 0$ by Corollary \ref{cor:alpha-invariant-injective}. Together with the vanishing of $H_\mathcal{G}^k(\nu_d^-, \nu_d')^a$, Proposition \ref{prop:vertical-LES} implies $\ker \zeta^k = 0$, and so $\ker \zeta^k = \ker \alpha^k = 0$ for $k > 4g-4-d_E+2d+1$.
Therefore, for all values of $k$ we have $\ker \alpha^k = \ker \zeta^k$.
\end{proof}

\section{Computation of the equivariant Betti numbers}\label{sec:computations}
Here we  use the results above, specifically Proposition \ref{prop:alpha}, together with the commutative diagram \eqref{eqn:diagram}, and derive an explicit formula for the equivariant Poincar\'e polynomial of $\mathcal{B}_0^{ss}(2,0)$ and $\mathcal{B}^{ss}(2,0)$.

We have the following  relationship between the equivariant Betti numbers of $X_d$ and $X_{d-1}$.
\begin{lem}\label{lem:zeta}
$$
\dim \, \ker \alpha^{k+1} - \dim \, \im \alpha^k = \dim H_{\mathcal{G}}^k(\nu_d', \nu_d'') - \dim H_{\mathcal{G}}^k(\nu_d^-, \nu_d'')
$$
\end{lem}

\begin{proof}

Using the vertical LES in diagram (\ref{eqn:diagram}) we have
\begin{align*}
\ker \zeta^{k+1} \cong \im \delta^k &\cong \frac{H_{\mathcal{G}}^k(\nu_d', \nu_d'')}{\ker \delta^k} \cong \frac{H_{\mathcal{G}}^k(\nu_d', \nu_d'')}{\im \lambda^k} \\
\im \lambda^k  &\cong \frac{H_{\mathcal{G}}^k(\nu_d^-, \nu_d'')}{\ker \lambda^k} \cong \frac{H_{\mathcal{G}}^k(\nu_d^-, \nu_d'')}{\im \zeta^k}
\end{align*}
Therefore
\begin{align*}
\dim \, \ker \zeta^{k+1} & = \dim H_{\mathcal{G}}^k(\nu_d', \nu_d'') - \dim \, \im \lambda^k \\
 & = \dim H_{\mathcal{G}}^k(\nu_d', \nu_d'') - \dim H_{\mathcal{G}}^k(\nu_d^-, \nu_d'') + \dim \, \im \zeta^k
\end{align*}
and so Proposition \ref{prop:alpha} implies
\begin{align*}
\dim \, \ker \alpha^{k+1} - \dim \, \im \alpha^k & = \dim \, \ker \zeta^{k+1} - \dim \, \im \zeta^k \\
 & = \dim H_{\mathcal{G}}^k(\nu_d', \nu_d'') - \dim H_{\mathcal{G}}^k(\nu_d^-, \nu_d'').
\end{align*}
completing the proof.
\end{proof}

\begin{prop}\label{prop:inductive-formulae}
$$
\dim H_{\mathcal{G}}^k(X_d) - \dim H_{\mathcal{G}}^k(X_{d-1}) = \dim H_{\mathcal{G}}^k(\nu_d^-, \nu_d'') - \dim H_{\mathcal{G}}^k(\nu_d', \nu_d'')
$$
In the fixed determinant case
\begin{multline}\label{eqn:fixed-determinant-betti-numbers}
\dim H_{\mathcal{G}}^k(X_d) - \dim H_{\mathcal{G}}^k(X_{d-1}) \\
 = \dim H^{k-2\mu_d} \left( J_d(M) \times BU(1)) \right) - \dim H^{k-2\mu_d}(\widetilde{S}^{2g- 2 + d_E - 2d} M).
\end{multline}
In the non-fixed determinant case 
\begin{multline}\label{eqn:non-fixed-betti-numbers}
\dim H_{\mathcal{G}}^k(X_d) - \dim H_{\mathcal{G}}^k(X_{d-1}) \\
 = \dim H^{k-2\mu_d}\left(J_d(M) \times J_n(M) \times BU(1) \times BU(1) \right) \\
 - \dim H^{k-2\mu_d} \left(S^{2g-2 + d_E - 2d} M \times J_d(M) \times BU(1) \right).
\end{multline}
\end{prop}

\begin{proof}
Lemma \ref{lem:zeta} shows that
\begin{align*}
\dim H_{\mathcal{G}}^k(X_d) &- \dim H_{\mathcal{G}}^k(X_{d-1}) \\
& = \dim \, \im \beta^k + \dim \, \ker \beta^k 
 - \dim \, \im \gamma^k - \dim \, \ker \gamma^k \\
 & = \dim \, \im \beta^k + \dim \, \im \alpha^k  - \dim \, \ker \alpha^{k+1} - \dim \, \im \beta^k \\
 & = \dim \, \im \alpha^k - \dim \, \ker \alpha^{k+1} \\
 & = \dim H_{\mathcal{G}}^k(\nu_d^-, \nu_d'') - \dim H_{\mathcal{G}}^k(\nu_d', \nu_d'') .
\end{align*}
In the fixed determinant case  use eqs.\  \eqref{eqn:nucohomology-doubleprime}, \eqref{eqn:nucohomology-prime-doubleprime}, \eqref{eqn:nonfixedcohomologycritical} and \eqref{eqn:nonfixednucohomology-prime-doubleprime} to obtain  \eqref{eqn:fixed-determinant-betti-numbers}. In the non-fixed determinant case use eqs.\   \eqref{eqn:nucohomology-doubleprime},  \eqref{eqn:nucohomology-prime-doubleprime},\eqref{eqn:fixedcohomologycritical} and \eqref{eqn:fixednucohomology-prime-doubleprime} to obtain \eqref{eqn:non-fixed-betti-numbers}.
\end{proof}

Inductively computing $H_{\mathcal{G}}^*(X_d)$ in terms of $H_{\mathcal{G}}^*(X_{d-1})$ for each value of $d$, we obtain the

\begin{proof}[Proof of Theorem \ref{thm:bettinumbers1}]

First we study the fixed determinant case. Eq. \eqref{eqn:fixed-determinant-betti-numbers} shows that  in both the degree zero and degree one case we have
\begin{equation*}\label{eqn:sum-inductive-fixed}
P_t^\mathcal{G} (\mathcal{B}) - P_t^\mathcal{G}(\mathcal{B}_0^{ss}) = \sum_{d = 1}^\infty t^{2 \mu_d} \frac{(1+t)^{2g}}{1-t^2} - \sum_{d=1}^{g-1} t^{2\mu_d} P_t (\widetilde{S}^{2g-2+d_E-2d} M)
\end{equation*}
 where $\mu_d = g-1+2d-d_E$. Note that the second sum has only $g-1$ terms because $H_\mathcal{G}^*(\nu_d', \nu_d'')$ is only non-zero if the vector space $H^0(L_1^* L_2 \otimes K)$ is non-zero, i.e. $d_E-2d+2g-2 \geq 0$, where $\deg L_1 = d$ and $\deg L_2 = d_E-d$.
 
Re-arranging this equation and substituting $P_t^\mathcal{G}(\mathcal{B}) = P_t(B \mathcal{G})$,
\begin{equation*}
P_t^\mathcal{G}(\mathcal{B}_0^{ss}) = P_t(B\mathcal{G}) - \sum_{d = 1}^\infty t^{2 \mu_d} \frac{(1+t)^{2g}}{1-t^2} + \sum_{d=1}^{g-1} t^{2\mu_d} P_t (\widetilde{S}^{2g-2+d_E-2d} M)
\end{equation*}
which proves \eqref{eqn:equivariantbetti-fixed0}. A similar argument using \eqref{eqn:non-fixed-betti-numbers} in Proposition \ref{prop:inductive-formulae} proves \eqref{eqn:equivariantbetti-nonfixed0}.
\end{proof}
 As mentioned in the Introduction, in the degree one case this gives a new proof of \cite[Thm.\ 7.6 (iv)]{Hitchin87} (fixed determinant case) and the results of \cite{HauselThaddeus04} (non-fixed determinant case).

In \cite[Sect.\  7]{Hitchin87} an explicit formula is given for the sum $$\sum_{d=1}^{g-1} t^{2\mu_d} P_t(\widetilde{S}^{2g-2d-1} M)$$ for $\mu_d = g+2d-2$, corresponding to the case where $\deg (E) = 1$. For the degree zero case we use eqs.\ \eqref{eqn:equivariantbetti-fixed0} and \eqref{eqn:equivariantbetti-nonfixed0}, together with the techniques of \cite{Hitchin87} to give the 

\begin{proof}[Proof of Corollary \ref{thm:maintheoremdegreezero}]

First, recall from \cite[Section 2]{AtiyahBott83} that for the rank $2$ fixed determinant case
\begin{equation}\label{eqn:atiyahbottc1=0}
P_t(B\mathcal{G}) = \frac{(1+t^3)^{2g}}{(1-t^2)(1-t^4)} ,
\end{equation}
and for the non-fixed determinant case
\begin{equation}\label{eqn:atiyahbottc1=0nonfixed}
P_t(B \mathcal{G}) = \frac{(1+t)^{2g} (1+t^3)^{2g}}{(1-t^2)^2 (1-t^4)} .
\end{equation}
Note that using the results from \cite[eq.\ (7.13)]{Hitchin87}, the last term in (\ref{eqn:equivariantbetti-fixed0}) is given by
\begin{align}\label{eqn:coversymmproduct}
\begin{split}
\sum_{d=1}^{g-1} t^{2\mu_d} P_t(\widetilde{S}^{2g-2d-2} M) & = \sum_{d=1}^{g-1} t^{2(g+2d-1)} P_t(S^{2g-2d-2} M) \\
 & \quad \quad + (2^{2g}-1) \sum_{d=1}^{g-1} \left( \begin{matrix} 2g-2 \\ 2g-2d-2 \end{matrix} \right) t^{4g+2d-4} \\
 & = \sum_{d=1}^{g-1} t^{2(g+2d-1)} P_t(S^{2g-2d-2} M) \\
 & \quad \quad + (2^{2g}-1) t^{4g-4} \sum_{d=1}^{g-1} \left( \begin{matrix} 2g-2 \\ 2g-2d-2 \end{matrix} \right) t^{2d}
\end{split}
\end{align}
Using the binomial theorem, the second term is
\begin{equation}\label{eqn:binomialterm}
\frac{1}{2}(2^{2g}-1) t^{4g-4} \left( (1+t)^{2g-2} + (1-t)^{2g-2} - 2 \right)
\end{equation}
The first term is calculated in the following lemma
\begin{lem} \label{eqn:evenbradlow}
\begin{align*}
\sum_{d=1}^{g-1} t^{2(g+2d-1)} P_t(S^{2g-2d-2} M) = &  -t^{4g-4} + \frac{t^{2g+2} (1+t)^{2g}}{(1-t^2)(1-t^4)} + \frac{(1-t)^{2g} t^{4g-4}}{4(1+t^2)} \\
& \negthickspace \negthickspace \negthickspace \negthickspace \negthickspace \negthickspace \negthickspace \negthickspace \negthickspace \negthickspace \negthickspace \negthickspace - \frac{(t+1)^{2g} t^{4g-4}}{2(t^2 - 1)} \left( \frac{2g}{t+1} + \frac{1}{t^2 - 1} - \frac{1}{2} + (3-2g) \right)
\end{align*}
\end{lem}

Part (b) of Corollary \ref{thm:maintheoremdegreezero} immediately follows from eqs.\ \eqref{eqn:equivariantbetti-nonfixed0}, \eqref{eqn:atiyahbottc1=0nonfixed} and \eqref{eqn:evenbradlow}. Part (a) follows from combining eqs.\ \eqref{eqn:equivariantbetti-fixed0}, \eqref{eqn:atiyahbottc1=0}, \eqref{eqn:coversymmproduct} and \eqref{eqn:binomialterm} and Lemma \ref{eqn:evenbradlow}.
\renewcommand{\qedsymbol}{}\end{proof}

\begin{proof}[Proof of Lemma \ref{eqn:evenbradlow}]
By \cite{MacDonald62}, $P_t(S^{2g-2d-2} M)$ is the coefficient of $x^{2g-2d-2}$ in $\frac{(1+xt)^{2g}}{(1-x)(1-xt^2)}$, or equivalently the coefficient of $x^{2g}$ in $\frac{x^{2d+2}(1+xt)^{2g}}{(1-x)(1-xt^2)}$. Therefore the sum 
\begin{equation*}
\sum_{d=1}^{g-1} t^{2(g+2d-1)} P_t(S^{2g-2d-2} M)
\end{equation*}
is the coefficient of $x^{2g}$ in
\begin{equation*}\label{eqn:finitesum}
\sum_{d=1}^{g-1} t^{2(g+2d-1)} x^{2d+2} \frac{(1+xt)^{2g}}{(1-x)(1-xt^2)}
\end{equation*}
which is equal to the coefficient of $x^{2g}$ in the following infinite sum
\begin{equation*}
\sum_{d=1}^\infty t^{2(g+2d-1)} x^{2d+2} \frac{(1+xt)^{2g}}{(1-x)(1-xt^2)}
\end{equation*}
The sum above is equal to
\begin{align*}
\sum_{d=1}^\infty t^{2(g+2d-1)} x^{2d+2} \frac{(1+xt)^{2g}}{(1-x)(1-xt^2)} & = t^{2g+2} x^4 \frac{(1+xt)^{2g}}{(1-x)(1-xt^2)} \sum_{d=1}^\infty (xt^2)^{2d-2} \\
 & = \frac{t^{2g+2} x^4 (1+xt)^{2g}}{(1-x)(1-xt^2)(1-x^2t^4)}
\end{align*}
Therefore the coefficient of $x^{2g}$ in the above sum is equal to the residue at $x=0$ of the function
\begin{equation*}
f(x) = \frac{(1+xt)^{2g} t^{2g+2}}{(1-x)(1-xt^2)^2 (1+xt^2)} \cdot \frac{1}{x^{2g-3}}
\end{equation*}
As in \cite{Hitchin87}, this residue can be computed in terms of the residues at the simple poles $x=1$ and $x=-t^{-2}$, the residue at the double pole $x=t^{-2}$, and the integral of $f(x)$ around a contour containing all of the poles. In this case the same methods can be used to compute the residues. However, unlike the situation in \cite{Hitchin87}, the contour integral is not asymptotically zero as the contour approaches the circle at infinity, so this must be computed here as well. To compute the integral, let $C_r$ be the circle of radius $r$ in the complex plane where $r>1$ and $r>t^{-2}$ (i.e.\ the disk inside $C_r$ contains all the poles of $f(x)$). Then for $\left| x \right| = r$ we have the following Laurent expansion of $f(x)$ centred at $x=0$.
\begin{align*}
\frac{(1+xt)^{2g} t^{2g+2} x^{3-2g}}{(1-x)(1-xt^2)^2(1+xt^2) } = & - \frac{(\frac{1}{x} + t)^{2g} t^{2g+2}}{xt^6 \left(1-\frac{1}{x} \right) \left(1- \frac{1}{xt^2} \right)^2 \left( 1 + \frac{1}{xt^2} \right) } \\
  = & -\frac{1}{xt^4} \left(\frac{t}{x} + t^2 \right)^{2g} \left( 1 + \frac{1}{x} + \cdots \right)  \\
 & \, \, \times \left( 1 + \frac{1}{xt^2} + \cdots \right)^2 \left( 1 - \frac{1}{xt^2} + \cdots \right) \\
 = & - \frac{t^{4g-4}}{x} + \mathrm{terms \, of \, order} \, x^{-n} \, \mathrm{where} \, n>1
\end{align*}
This series expansion is uniformly convergent on the annulus $\{ x : r-\varepsilon < x < r+\varepsilon \}$ for $r>1$, $r> t^{-2}$ and $\varepsilon$ small enough so that the closure of the annulus doesn't contain any of the poles of $f(x)$. As $r \rightarrow \infty$ the series asymptotically approaches $-t^{4g-4}/x$, and so the integral approaches
\begin{equation}\label{eqn:integral}
\lim_{r \rightarrow \infty} \frac{1}{2\pi i} \int_{C_r} \frac{(1+xt)^{2g} t^{2g+2} x^{3-2g}}{(1-x)(1-xt^2)^2(1+xt^2) } \, dx = -t^{4g-4}
\end{equation}
The residues of $f(x)$ at $x=1$, $x=-t^{-2}$ and $x=t^{-2}$ are similar to the results obtained in \cite{Hitchin87}. At the simple pole $x=1$,
\begin{equation}\label{eqn:resx=1}
Res_{x=1} f(x) = -\frac{t^{2g+2} (1+t)^{2g}}{(1-t^2)(1-t^4)}
\end{equation}
At the simple pole $x=-t^{-2}$
\begin{equation}\label{eqn:resx=-t2}
Res_{x=-t^{-2}} f(x) = - \frac{(1-t)^{2g} t^{4g-4}}{4(1+t^2)}
\end{equation}
and at the double pole $x=t^{-2}$
\begin{equation}\label{eqn:resx=t2}
Res_{x=t^{-2}} f(x) = \frac{(t+1)^{2g} t^{4g-4}}{2(t^2 - 1)} \left( \frac{2g}{t+1} + \frac{1}{t^2 - 1} - \frac{1}{2} + (3-2g) \right)
\end{equation}
Combining (\ref{eqn:integral}), (\ref{eqn:resx=1}), (\ref{eqn:resx=-t2}) and (\ref{eqn:resx=t2}) we have
\begin{align*}
\sum_{d=1}^{g-1} t^{2(g+2d-1)} & P_t(S^{2g-2d-2} M) = -t^{4g-4} + \frac{t^{2g+2} (1+t)^{2g}}{(1-t^2)(1-t^4)} + \frac{(1-t)^{2g} t^{4g-4}}{4(1+t^2)} \\
 & \, \, \, \, - \frac{(t+1)^{2g} t^{4g-4}}{2(t^2 - 1)} \left( \frac{2g}{t+1} + \frac{1}{t^2 - 1} - \frac{1}{2} + (3-2g) \right)
\end{align*}
thus completing the proof of the lemma and therefore also of Corollary \ref{thm:maintheoremdegreezero}.
\end{proof}

\bibliographystyle{plain}
\bibliography{ref}

\end{document}